\documentclass{amsart}
\usepackage[dvips]{graphicx}
\usepackage{amscd}
\usepackage{amsmath}
\usepackage{amsxtra}
\usepackage{amsfonts}
\usepackage{amssymb}

\newtheorem{theorem}{Theorem}[section]

\newtheorem{lemma}[theorem]{Lemma}
\newtheorem{proposition}[theorem]{Proposition}

\theoremstyle{definition}
\newtheorem{definition}[theorem]{Definition}
\newtheorem{remark}[theorem]{Remark}

\newtheorem*{question}{Question}
\newtheorem{example}[theorem]{Example}
\theoremstyle{remark}

\renewcommand{\theclaim}{\textup{\theclaim}}

\newtheorem*{acknowledgements}{Acknowledgements}

\numberwithin{equation}{section}

\def\openone

{\mathchoice

{\hbox{\upshape \small1\kern-3.3pt\normalsize1}}

{\hbox{\upshape \small1\kern-3.3pt\normalsize1}}

{\hbox{\upshape \tiny1\kern-2.3pt\SMALL1}}

{\hbox{\upshape \Tiny1\kern-2pt\tiny1}}}

\makeatletter

\newbox\ipbox

\newcommand{\ip}[2]{\left\langle #1\, , \,#2\right\rangle}
\newcommand{\diracb}[1]{\left\langle #1\mathrel{\mathchoice

{\setbox\ipbox=\hbox{$\displaystyle \left\langle\mathstrut
#1\right.$}

\vrule height\ht\ipbox width0.25pt depth\dp\ipbox}

{\setbox\ipbox=\hbox{$\textstyle \left\langle\mathstrut
#1\right.$}

\vrule height\ht\ipbox width0.25pt depth\dp\ipbox}

{\setbox\ipbox=\hbox{$\scriptstyle \left\langle\mathstrut
#1\right.$}

\vrule height\ht\ipbox width0.25pt depth\dp\ipbox}

{\setbox\ipbox=\hbox{$\scriptscriptstyle \left\langle\mathstrut
#1\right.$}

\vrule height\ht\ipbox width0.25pt depth\dp\ipbox}

}\right. }

\newcommand{\dirack}[1]{\left. \mathrel{\mathchoice

{\setbox\ipbox=\hbox{$\displaystyle \left.\mathstrut
#1\right\rangle$}

\vrule height\ht\ipbox width0.25pt depth\dp\ipbox}

{\setbox\ipbox=\hbox{$\textstyle \left.\mathstrut
#1\right\rangle$}

\vrule height\ht\ipbox width0.25pt depth\dp\ipbox}

{\setbox\ipbox=\hbox{$\scriptstyle \left.\mathstrut
#1\right\rangle$}

\vrule height\ht\ipbox width0.25pt depth\dp\ipbox}

{\setbox\ipbox=\hbox{$\scriptscriptstyle \left.\mathstrut
#1\right\rangle$}

\vrule height\ht\ipbox width0.25pt depth\dp\ipbox}

} #1\right\rangle}

\newcommand{\cj}[1]{\overline{#1}}

\newcommand{\bz}{\mathbb{Z}}

\newcommand{\br}{\mathbb{R}}
\newcommand{\bc}{\mathbb{C}}

\newcommand{\bn}{\mathbb{N}}

\def\blfootnote{\xdef\@thefnmark{}\@footnotetext}

\renewcommand{\mod}{\operatorname{mod}}

\input xy
\xyoption{all}
\usepackage{amssymb}

\newcommand{\Per}{\operatorname*{Per}}
\newcommand{\Span}{\overline{\mbox{span}}}

\def\-{^{-1}}

\begin{document}
\title[Orthonormal dilations of Parseval wavelets]{Orthonormal dilations of Parseval wavelets}
\author{Dorin Ervin Dutkay${}^*$}

\blfootnote{${}^*$ Research supported in part by a grant from the National Science Foundation DMS-0704191}
\address{[Dorin Ervin Dutkay] University of Central Florida\\
	Department of Mathematics\\
	4000 Central Florida Blvd.\\
	P.O. Box 161364\\
	Orlando, FL 32816-1364\\
U.S.A.\\} \email{ddutkay@mail.ucf.edu}
\author{Deguang Han}
\address{[Deguang Han] University of Central Florida\\
	Department of Mathematics\\
	4000 Central Florida Blvd.\\
	P.O. Box 161364\\
	Orlando, FL 32816-1364\\
U.S.A.\\}
\email{dhan@pegasus.cc.ucf.edu}
\author{Gabriel Picioroaga}
\address{[Gabriel Picioroaga] Department of Mathematical Sciences\\
Binghamton University\\
Binghamton, New York 13902-6000 \\
U.S.A.}
\email{gabriel@math.binghamton.edu}
\author{Qiyu Sun}
\address{[Qiyu Sun] University of Central Florida\\
	Department of Mathematics\\
	4000 Central Florida Blvd.\\
	P.O. Box 161364\\
	Orlando, FL 32816-1364\\
U.S.A.\\}
\email{qsun@mail.ucf.edu}
\thanks{} 
\begin{abstract}
We prove that any Parseval wavelet frame is the projection of an orthonormal wavelet basis for a representation of the Baumslag-Solitar group
$$BS(1,2)=\langle u,t\,|\, utu^{-1}=t^2\rangle.$$
We give a precise description of this representation in some special cases, and show that for wavelet sets, it is related to symbolic dynamics (Theorem \ref{thdilmra}). We prove that the structure of the representation depends on the analysis of certain finite orbits for the associated symbolic dynamics (Theorem \ref{th3_24}). We give concrete examples of Parseval wavelets for which we compute the orthonormal dilations in detail; we construct Parseval wavelet sets which have infinitely many non-isomorphic orthonormal dilations.
\end{abstract}
\maketitle \tableofcontents
\section{Introduction}
An {\it orthonormal wavelet} in $L^2(\br)$ is a function $\psi\in L^2(\br)$ with the property that
\begin{equation}\label{eqwv1}
\{2^{j/2}\psi(2^j\cdot-k)\,|\,j\in\bz, k\in\bz\}
\end{equation}
is an orthonormal basis for $L^2(\br)$. 

The theory of wavelets has found numerous applications in a variety of areas such as signal processing, image compression and numerical analysis \cite{Dau92,Mal98,BrJo02}. The geometry of orthonormal wavelets is fairly well understood \cite{MR1987107,HST07,MR1666001,MR1968123,MR1468204,MR1658652}. The main technique used in the study of orthonormal wavelets is the multiresolution analysis (MRA) introduced by Mallat and Meyer \cite{Dau92,BrJo02} and its generalizations \cite{BMM99,BCM02,BJMP05}.

A function $\psi$ is called a {\it Parseval wavelet} if the family in \eqref{eqwv1} is a Parseval frame for $L^2(\br)$. A {\it Parseval frame} for a Hilbert space $\mathcal H$ is a family $\{e_i\,|\, i\in I\}$ of vectors in $\mathcal H$ that satisfies the Parseval identity
$$\|f\|^2=\sum_{i\in I}|\ip{f}{e_i}|^2,\quad (f\in \mathcal H).$$
Parseval wavelets are more flexible than their orthogonal counterparts. They can have a certain degree of symmetry, which is advantageous in applications \cite{Dau92}. The redundancy in the associated basis decompositions can be useful in compression problems. However, since the Parseval wavelets are not orthogonal, they can have complicated correlations, and the multiresolution techniques cannot be applied. In this paper we will prove that Parseval wavelets can be obtained from orthonormal wavelets by Hilbert space dilations and projections, and therefore one can use multiresolutions in this case too.

There is a very general result linking orthonormal bases to Parseval frames \cite{HaLa00}, which says that every Parseval frame is a projection of an orthonormal basis. More precisely, if $\{e_i\,|\,i\in I\}$ is a Parseval frame for a Hilbert space $\mathcal H$, then there exists a bigger Hilbert space $\tilde{\mathcal H}\supset\mathcal H$ and an orthonormal basis $\{\tilde e_i\,|\,i\in I\}$ for $\tilde{\mathcal H}$ such that 
$P_{\mathcal H}\tilde e_i=e_i$ for all $i\in I$, where $P_{\mathcal H}$ is the orthogonal projection onto the subspace $\mathcal H$. We say that the Parseval frame $\{e_i\,|\,i\in I\}$ can be {\it dilated} to the orthonormal basis $\{\tilde e_i\,|\,i\in I\}$, or that $\{\tilde e_i\,|\,i\in I\}$ is an {\it orthonormal dilation} of the Parseval frame $\{e_i\,|\, i\in I\}$. We call $\{(1_{\tilde{\mathcal H}}-P_{\mathcal H})\tilde e_i\,|\,i\in I\}$ a {\it complementary frame} for $\{e_i\,|\, i\in I\}$.

Then a natural question is: if the Parseval frame $\{e_i\,|\,i\in I\}$ has some additional structure can we dilate it to an orthonormal basis that shares similar properties? In the case of frames generated by actions groups or for Gabor frames, the answer is positive \cite{HaLa00}. For Parseval wavelets there are some dilation results in the literature \cite{HaLa00,GuHa05,BDP05} which apply to some particular classes of wavelets. In this paper we give a complete solution for the general case and prove that the affine structure attached to the wavelet basis can be preserved under orthonormal dilations (Theorem \ref{thdilation}).

To formulate the question more explicitly, let us express the family in \eqref{eqwv1} in terms of the action of unitary operators. In $L^2(\br)$ we consider two unitary operators: the translation operator $T_0$, and the dilation operator $$ T_0f=f(\cdot-1),\quad U_0f=\frac{1}{\sqrt{2}}f\left(\frac{\cdot}{2}\right),\quad (f\in L^2(\br)).$$ Then the family in \eqref{eqwv1} is 
\begin{equation}\label{eqwv2}
\{U_0^jT_0^k\psi\,|\,j,k\in\bz\}.
\end{equation}
The two operators $U_0$ and $T_0$ satisfy the relation $U_0T_0U_0^{-1}=T_0^2$, therefore we are dealing with a unitary representation of the Baumslag-Solitar group with two generators and one relation:
\begin{equation}\label{eqbs}
BS(1,2):=\langle u,t\,|\,utu^{-1}=t^2\rangle.
\end{equation}
Any representation of the Baumslag-Solitar group $BS(1,2)$ is in fact given by two unitary operators $U$ and $T$ on some Hilbert space $\mathcal H$, that satisfy the relation $UTU^{-1}=T^2$.

While the Baumslag-Solitar group appears to be quite simple, this can be deceiving, there are several extremely interesting results about it in the literature which reveal surprising properties. The Baumslag-Solitar group is of independent interest in combinatorial topology and operator algebras. In \cite{MR1780213} the authors compute the spectrum of the Markov operator associated to this group, basing their result on the Generalized Riemann Hypothesis! In \cite{MR1608595,MR1709862} these groups are shown to satisfy some rigidity properties, and at the same time they are not lattices in Lie groups.

\begin{definition}\label{defwv}
Let $\{U,T\}$ be a representation of the Baumslag-Solitar group $BS(1,2)$ on a Hilbert space $\mathcal H$. We call a vector $\psi\in\mathcal H$ a {\it Parseval (orthonormal) wavelet} for this representation, if 
$$\{U^jT^k\psi\,|\,j,k\in\bz\}$$
is a Parseval frame (orthonormal basis) for $\mathcal H$.
\end{definition}

 To distinguish between the two levels, the initial problem, and the
extended version for the dilated Hilbert space, we use the name
``super-representation'' for the latter. In the Parseval-wavelet case, the
dilated version acquires the orthonormal structure, while still preserving
the affine scaling relations dictated by the Baumslag-Solitar group.

\begin{definition}
Let $\{\tilde U,\tilde T\}$, $\{U,T\}$ be representations of the group $BS(1,2)$ on the Hilbert spaces $\tilde{\mathcal H}$ and $\mathcal H$ respectively. We say that $\{U,T\}$ is a subrepresentation of $\{\tilde U,\tilde T\}$ if $\mathcal H\subset\tilde{\mathcal H}$, the projection $P_{\mathcal H}$ onto the subspace $\mathcal H$ commutes with $\tilde U$ and $\tilde T$ and $P_{\mathcal H}\tilde UP_{\mathcal H}=U$, $P_{\mathcal H}\tilde TP_{\mathcal H}=T$. We will also say that $\{\tilde U,\tilde T\}$ is a {\it super-representation} of $\{U,T\}$.
\end{definition}

Now we can formulate our question more precisely:
\begin{question}{\bf (Dilations of Parseval wavelets)}
{\it Let $\psi$ be a Parseval wavelet for a representation $\{U,T\}$ of the group $BS(1,2)$ on a Hilbert space $\mathcal H$. Does there exist a representation $\{\tilde U,\tilde T\}$ of $BS(1,2)$ on a bigger space $\tilde{\mathcal H}\supset\mathcal H$ and a vector $\tilde\psi\in\tilde{\mathcal H}$ such that $\{U,T\}$ is a subrepresentation of $\{\tilde U,\tilde T\}$, $\tilde \psi$ is an orthonormal wavelet for $\{\tilde U, \tilde T\}$, and $P_{\mathcal H}\tilde\psi=\psi$? $P_{\mathcal H}$ is the orthogonal projection onto the subspace $\mathcal H$.}
\end{question}

We will give a positive answer to our Question in Theorem \ref{thdilation}: any Parseval wavelet can be dilated to an orthonormal wavelet. 

The results from \cite{HaLa00,Du04} cannot be applied directly because the family in \eqref{eqwv2} does not involve the entire group $BS(1,2)$, but only a subset of it, namely the elements of the form $u^jt^k$ with $j,k\in\bz$. Our construction of orthonormal dilations will be based on the general theory of Hilbert spaces built out of positive definite functions. This is essentially contained in Theorem \ref{thpdtobs}. Since seminal papers by Krein and Rudin, the problem of finding a positive extension for a positive definite map from a subset to the entire group is known to be notoriously difficult, few results are available, each for a very particular case, see \cite{MR0004333,MR0151796,MR931979,MR1794293, Jor89,Jor90, Jor91}. Theorem \ref{thpdtobs} adds one to the list: a positive definite map on the subset of $BS(1,2)$ determined by the wavelet family in \eqref{eqwv2} can be extended to the whole group $BS(1,2)$. 

\medskip

There are several abstract precursors to our extension problems. These include: unitary dilation of isometries, Stinespring dilations in operator algebras, or Naimark dilations for operator valued measures. However, these earlier results lack computational detail. Our results in Section \ref{waveletsets} identify the right abstract models and provide algorithms for the computation of the orthogonal dilation. 

\begin{question}{\bf (Explicit dilations)}
{\it Let $\psi$ be a Parseval wavelet for $\{U_0,T_0\}$ in $L^2(\br)$ and let $\tilde\psi$ be an orthonormal wavelet for a super-representation
$\{\tilde U,\tilde T\}$. What is the precise structure of the representation $\{\tilde U,\tilde T\}$ and what is $\tilde\psi$?}
\end{question}

We do not have a complete answer to this question; nevertheless, we are be able to construct a concrete orthonormal dilation in the special case of Parseval wavelet sets, which has the advantage that preserves the multiresolution structure. We believe that our results can be extended to more complicated Parseval wavelets. We offer a computational correspondence
between two seemingly unrelated areas, representations of the
Baumslag-Solitar group on the one hand, and on the other a formula
for the geometry and for invariants of wavelet sets (sections \ref{waveletsets} and \ref{ex}).

We recall some of the concrete dilation results in the literature. 
A {\it wavelet set} is a wavelet $\psi$ such that its Fourier transform $\hat\psi$ is a characteristic function. In \cite{HaLa00,GuHa05} many examples of Parseval wavelet sets are provided where the orthonormal dilation lies in the space $L^2(\br)\oplus\dots \oplus L^2(\br)$ with the representation of the group $BS(1,2)$ given by a simple amplification of the representation $\{U_0,T_0\}$ in $L^2(\br)$:
$\tilde U=U_0\oplus \dots \oplus U_0$, $\tilde T=T_0\oplus\dots \oplus T_0$. There is one issue with this representation, as shown in \cite{HaLa00}: it does not have orthogonal multiresolution wavelets. Therefore if we start from a multiresolution Parseval wavelet in $L^2(\br)$, and we want to dilate it to an orthonormal wavelet in such a way that this super-wavelet comes also from a multiresolution, then we have to look somewhere else, and replace the amplification with another representation.

An answer to this problem can be found in \cite{BDP05}. We illustrate it by a classical example: the stretched Haar wavelet
$\psi=\frac12\chi_{[0,3/2)}-\frac12\chi_{[3/2,3)}$ is a non-orthogonal Parseval wavelet that is constructed from a multiresolution with low-pass filter $m_0(x)=\frac{1}{\sqrt2}(1+e^{2\pi i 3x})$ and scaling function $\varphi=\frac13\chi_{[0,3)}$ (see \cite{Dau92,BrJo02}). In \cite{BDP05} it was shown that, in order to construct the orthonormal dilation wavelet that preserves the multiresolution, one has to consider the representation of $BS(1,2)$ on $L^2(\br)\oplus L^2(\br)\oplus L^2(\br)$ given by
$$T_3(f_1,f_2,f_3)=(T_0f_1,e^{2\pi i/3}T_0f_2,e^{4\pi i/3}T_0f_3),\quad
U_3(f_1,f_2,f_3)=(U_0f_1,U_0f_3,U_0f_2).$$
The dilated orthonormal wavelet is
$$\tilde\psi=(\psi,\psi,\psi).$$
This is a multiresolution wavelet that has the associated scaling function $\tilde\varphi=(\varphi,\varphi,\varphi)$.

The theory in \cite{BDP05} shows that this procedure works in a more general case, e.g. when $\psi$ is a compactly supported multiresolution Parseval frame. Then the orthonormal dilation can be realized in a similar ``twisted'' amplification of the representation of the group $BS(1,2)$ in $L^2(\br)$.

One difficulty of the theory in \cite{BDP05} is that it requires the low-pass filter to have a finite number of zeros, and therefore, it cannot be used for Parseval wavelet sets. In section \ref{waveletsets} we will construct orthonormal dilations in the very special case of multiresolution Parseval wavelet sets, and show that even in this particular case there are interesting connections to symbolic dynamics. We will show that the orthonormal dilations are by no means unique, and in Proposition \ref{propsemi} we prove that, when the Parseval wavelet is semi-orthogonal, the dilation can be realized in a subrepresentation of $L^2(\br)\oplus L^2(\br)$ with $U=U_0\oplus U_0$ and $T=T_0\oplus T_0$.

Here is the summary of the paper: in Section \ref{general} we analyze Parseval wavelets for abstract representations of the group $BS(1,2)$. We show that every Parseval wavelet can be dilated to an orthogonal wavelet (Theorem \ref{thdilation}). This is based on the fact that positive definite maps on the subset of the group $BS(1,2)$ can be extended to a positive definite map on the entire group (Theorem \ref{thpdtobs}). We prove that if the Parseval wavelet is semi-orthogonal, then we have a concrete form of this orthonormal dilation, explicitly, the Parseval wavelet has a complementary wavelet in a subspace of $L^2(\br)$ (Proposition \ref{propsemi}).

In Section \ref{waveletsets} we shift our focus to MRA Parseval wavelet sets in $L^2(\br)$ and we give a concrete form of an orthonormal dilation (Theorem \ref{thdilmra}). This requires several steps: in Section \ref{scaling} we show how a low-pass filter and scaling function can be constructed for a MRA Parseval wavelet set. The low-pass filter is then used in Section \ref{dilated} to construct a representation of the $BS(1,2)$ on a symbolic space. This representation  will contain the orthonormal dilation. In Section \ref{encoding} we show how the classical representation on $L^2(\br)$ can be embedded in this representation. Section \ref{main} contains our dilation result for Parseval wavelet sets. Theorem \ref{thdilmra} provides the concrete orthonormal dilation for a Parseval wavelet sets which preserves the multiresolution structure. In Section \ref{cyclic} we show that, under certain assumptions on the low-pass filter, the dilated representation of Section \ref{dilated} is in fact the same as the one used in \cite{BDP05}, of the type we mentioned above for the stretched Haar wavelet. The representation is based on cyclic orbits for the associated symbolic dynamics. 

In the final section of the paper we give some concrete examples of orthonormal dilations of Parseval wavelet sets. In Example \ref{ex1} we show that the family of Parseval wavelet sets $\hat\psi_{[-2a,-a]\cup[a,2a]}$, with $0<a\leq\frac14$ can be dilated in the same representation of $BS(1,2)$ as the stretched Haar wavelet. In Example \ref{ex2} we construct an orthonormal dilation of $\hat\psi_{[-\frac14,-\frac18]\cup[\frac18,\frac14]}$ in a different representation, thus proving that the orthonormal dilation is not unique. Example \ref{ex3} proves that, in some cases, the cycles are not sufficient to describe the orthonormal dilation, therefore the results of Section \ref{cyclic} do not give a complete picture of the possible representations of Theorem \ref{thdilmra}. In Example \ref{ex4} we prove that if $a$ is small enough, the Parseval wavelet set in Example \ref{ex1} has infinitely many non-isomorphic orthonormal dilations.

\section{General dilations of Parseval wavelets}\label{general}
We want to construct an orthonormal dilation of a Parseval wavelet. For this we will first construct a certain positive definite map, following the ideas in \cite{Du04}. Recall that a map
$K:X\times X\rightarrow\bc$ is said to be {\it positive definite} if for all finite sets $F\subset X$ and any $x_i\in X, \xi_i\in\bc$, with $i\in F$ one has
$$\sum_{i,j\in F}K(x_i,x_j)\xi_i\cj\xi_j\geq0.$$
From the positive definite map $K$, one can construct a Hilbert space and a family of vectors that have the inner products determined by $K$. Then the crucial point is to construct the unitary operators $U$ and $T$ and a $\psi$ such that this family of vectors is equal to $\{U^jT^k\psi\,|\, j,k\in\bz\}$ and such that the relation $UTU^{-1}=T^2$ is satisfied.

\begin{theorem}\label{thpdtobs}
Let $K:\bz^2\times\bz^2\rightarrow\bc$ be positive definite, and assume that the following conditions are satisfied:
\begin{equation}\label{eqinv1}
K((j,k),(j',k'))=K((j+1,k),(j'+1,k')),\quad(j,j',k,k'\in\bz),
\end{equation}
and
\begin{equation}\label{eqinv2}
K((j,k),(j',k'))=K((j,2^{-j}+k),(j',2^{-j'}+k')),\quad(j,j'\leq0, k,k'\in\bz).
\end{equation}

Then there exists a Hilbert space $H$, a representation $U$, $T$ of the Baumslag-Solitar group $BS(1,2)$, and a vector $\psi\in H$ such that
$$\ip{U^jT^k\psi}{U^{j'}T^{k'}\psi}=K((j,k),(j',k')),\quad(j,j',k,k'\in\bz)$$
and $\overline{\mbox{span}}\{U^jT^k\psi\,|\,j,k\in\bz\}=H$.
\end{theorem}
\begin{remark}
Before we give the proof of this theorem, let us explain where the relations \eqref{eqinv1} and \eqref{eqinv2} come from. If $U, T$ is a representation of the Baumslag-Solitar group $BS(1,2)$ on a Hilbert space $H$ and $\psi\in H$, then a simple computation that uses the fact that $U$ and $T$ are unitary and $TU^j=U^jT^{2^{-j}}$ for $j\leq0$, shows that the map
$$K((j,k),(j',k'))=\ip{U^jT^k\psi}{U^{j'}T^{k'}\psi},\quad(j,j',k,k'\in\bz)$$
satisfies \eqref{eqinv1} and \eqref{eqinv2}, and it is of course positive definite.
\end{remark}

\begin{proof}
Using Kolmogorov's result mentioned in \cite[Theorem 2.2]{Du04}, we obtain a Hilbert space $H$ and a map $v:\bz^2\rightarrow H$ such that $\ip{v(j,k)}{v(j',k')}=K((j,k),(j',k'))$ for all $j,j',k,k'$, and such that $\overline{\mbox{span}}\{v(j,k)\,|\,j,k\in\bz\}=H$. 
\par
Define the operator $U$ on $H$, by $Uv(j,k)=v(j+1,k)$ for all $j,k\in\bz$, and extend linearly. Then we claim that $U$ is a unitary operator. 
\par
Indeed, we have for all finite subsets $F$ of $\bz^2$ and $\alpha_{j,k}\in\bc$:
$$\left\|U(\sum_{(j,k)\in F}\alpha_{j,k}v(j,k))\right\|^2=\sum_{(j,k),(j',k')\in F}\alpha_{j,k}\cj{\alpha_{j',k'}}K((j+1,k),(j'+1,k'))=(\ast)$$
and using \eqref{eqinv1},
$$(\ast)=\sum_{(j,k),(j',k')\in F}\alpha_{j,k}\cj{\alpha_{j',k'}}K((j,k),(j',k'))=\left\|\sum_{(j,k)\in F}\alpha_{j,k}v(j,k)\right\|^2.$$
This shows that $U$ is an isometry, and since $\Span\{v(j+1,k)\,|\,j,k\in\bz\}=H$, $U$ is unitary. 
\par
Let $\tilde V_l:=\Span\{v(j,k)\,|\,j\leq l,k\in\bz\}$ for all $l\geq -1$. Then $\tilde V_l\subset\tilde V_{l+1}$, and $U\tilde V_l=\tilde V_{l+1}$, for $l\geq-1$. Note also that $\overline{\cup_l\tilde V_l}=H$. 
\par
Define the operator $T_0$ on $\tilde V_0$ by $T_0v(j,k)=v(j,2^{-j}+k)$, for $j\leq0,k\in\bz$. We check that $T_0$ is an isometry on $\tilde V_0$. 
$$\left\|\sum_{(j,k)\in F}\alpha_{j,k}v(j,2^{-j}+k)\right\|^2=\sum_{(j,k),(j',k')\in F}\alpha_{j,k}\cj{\alpha_{j',k'}}K((j,2^{-j}+k),(j',2^{-j'}+k))=(\ast)$$
and using \eqref{eqinv2},
$$(\ast)=\sum_{(j,k),(j',k')\in F}\alpha_{j,k}\cj{\alpha_{j',k'}}K((j,k),(j',k'))=\left\|\sum_{(j,k)\in F}\alpha_{j,k}v(j,k)\right\|^2.$$
This proves that $T_0$ is an isometry. 
\par
Clearly $T_0\tilde V_0=\tilde V_0$ and $T_0\tilde V_{-1}=\tilde V_{-1}$. 
\par
Define $\tilde W_l:=\tilde V_l\ominus\tilde V_{l-1}$ for $l\geq 0$. Then $H=\tilde V_0\oplus\oplus_{l\geq 1}\tilde W_l$. Also, since $T_0$ is unitary on $\tilde V_0$, $T_0\tilde W_0=\tilde W_0$, and, since $U$ is unitary, $U\tilde W_l=\tilde W_{l+1}$ for all $l\geq 0$. 
\par
We will need the following lemma, which can be easily obtained by an application of Borel functional calculus:
\begin{lemma}\label{lemsquni}
If $a$ is a unitary operator on a Hilbert space then there exists a unitary operator $b$, on the same Hilbert space, such that $b^2=a$.
\end{lemma}
Now define $T_1:\tilde W_1\rightarrow\tilde W_1$ as follows: the operator $UT_0U^{-1}$ is unitary on $\tilde W_1$, so by Lemma \ref{lemsquni}, there exists a unitary operator $T_1$ on $\tilde W_1$ such that $T_1^2=UT_0U^{-1}$. 
\par
By induction, we use Lemma \ref{lemsquni} to define the unitary operator $T_l$ on $\tilde W_l$ such that $T_l^2=UT_{l-1}U^{-1}$. 
\par
Now we can define the unitary operator $T$ on $H$ such that $T$ on $\tilde V_0$ is $T_0$, and $T$ on $\tilde W_l$ is $T_l$ for all $l\geq 1$.
\par
We check that $UTU^{-1}=T^2$. First, on $\tilde V_0$: take $j\leq0,k\in\bz$.
$$UTU^{-1}v(j,k)=UT_0v(j-1,k)=Uv(j-1,2^{-j+1}+k)=v(j,2^{-j+1}+k),$$
$$T^2v(j,k)=T_0T_0v(j,k)=T_0v(j,2^{-j}+k)=v(j,2^{-j}+2^{-j}+k)=v(j,2^{-j+1}+k).$$
Then, on $\tilde W_l$ for $l\geq 1$, $UTU^{-1}=UT_{l-1}U^{-1}=T_l^2=T$. 
\par
Let $\psi:=v(0,0)$. Then $U^jT^k\psi=U^jv(0,k)=v(j,k)$ and everything follows.
\end{proof}

To construct the orthonormal dilation of the Parseval wavelet and the associated representation of the Baumslag-Solitar group, we will find what the positive definite map associated to the ``complement'' should be. If $\{e_i\,|\,i\in I\}$ is a Parseval frame and $\{\tilde e_i\,|\,i\in I\}$ is an orthonormal dilation, then the complement is $\{\tilde e_i-e_i\,|\,i\in I\}$, so the positive definite map associated to the complement is $K_2(i,j)=\ip{\tilde e_i-e_i}{\tilde e_j-e_j}=\delta_{i,j}-\ip{e_i}{e_j}$, for $i,j\in I$. (We used here the fact that $\tilde e_i-e_i$ is orthogonal to $e_j$ for all $i,j\in I$.)

Using the positive definite map of the complement, we can construct a ``complementary'' representation of the group $BS(1,2)$ and the complementary Parseval frame $\psi_2$. The orthonormal dilation is then obtained by a direct sum of the two components. The details of these steps are contained in the following lemmas.

\begin{lemma}\label{lemma2}
If $(f_i)_{i\in I}$ is a Parseval frame for a Hilbert space $H$ then, for all $F\subset I$ finite and all $\alpha_i\in\bc$, $(i\in F)$, 
$$\left\|\sum_{i\in F}\alpha_if_i\right\|^2\leq \sum_{i\in F}|\alpha_i|^2.$$
This implies that 
$$K(i,j):=\delta_{i,j}-\ip{f_i}{f_j},\quad(i,j\in I)$$
is positive definite. 
\end{lemma}

\begin{proof}
By \cite{HaLa00} there exist a Hilbert space $K$, $K\supset H$, and $(e_i)_{i\in I}$ an orthonormal basis for $K$, such that, if $P$ is the projection onto $H$, then $Pe_i=f_i$ for all $i\in I$. Then
$$\left\|\sum_{i\in F}\alpha_if_i\right\|^2=\left\|P(\sum_{i\in F}\alpha_ie_i)\right\|^2\leq \left\|\sum_{i\in F}\alpha_ie_i\right\|^2=\sum_{i\in F}|\alpha_i|^2.$$
We check that $K$ is positive definite:
$$\sum_{i,j\in F}\alpha_i\cj\alpha_jK(i,j)=\sum_{i,j\in F}(\delta_{i,j}-\ip{f_i}{f_j})\alpha_i\cj\alpha_j=\sum_{i\in F}|\alpha_i|^2-\left\|\sum_{i\in F}\alpha_if_i\right\|^2\geq0.$$
\end{proof}

\begin{lemma}\label{lemma3}
If $(f_i)_{i\in I}$ is a Parseval frame for $H_1$ and $(g_i)_{i\in I}$ are vectors that span $H_2$ such that $\ip{g_i}{g_j}=\delta_{i,j}-\ip{f_i}{f_j}$ for all $i,j\in I$, then $(f_i\oplus g_i)_{i\in I}$ is an orthonormal basis for $H_1\oplus H_2$, and $(g_i)_{i\in I}$ is a Parseval frame for $H_2$.
\end{lemma}

\begin{proof}
As in \cite[Corollary 1.3]{HaLa00}, consider a strong complementary Parseval frame $(\tilde g_i)_{i\in I}$, i.e., $(f_i\oplus\tilde g_i)_{i\in I}$ is an orthonormal basis for $H_1\oplus\tilde H_2$, for some Hilbert space $\tilde H_2$. Then $\ip{\tilde g_i}{\tilde g_j}=\delta_{i,j}-\ip{f_i}{f_j}=\ip{g_i}{g_j}$. Define the operator $W$ from $\tilde H_2$ to $H_2$, by $W\tilde g_i=g_i$. Clearly, $W$ is an isometry and, since $\Span\{g_i\,|\,i\in I\}=H_2$, it follows that $W$ is unitary, so $(g_i)_{i\in I}$ is also a Parseval frame. 
\par
The operator $I\oplus U$ is unitary so $(f_i\oplus g_i)_{i\in I}$ is an orthonormal basis for $H_1\oplus H_2$.  
\end{proof}

With Lemma \ref{lemma2} and Lemma \ref{lemma3} the desired dilation result follows:

\begin{theorem}\label{thdilation}
Any Parseval wavelet can be dilated to an orthonormal wavelet. More precisely, let $\{U, T\}$ be a representation of the Baumslag-Solitar group $BS(1,2)$ on some Hilbert space $H$. Let $\psi$ be a Parseval wavelet for $\{U,T\}$ on $H$. Then there exists a Hilbert space $H_2$, a representation $\{U_2$, $T_2\}$ of the group $BS(1,2)$, and a Parseval wavelet $\psi_2$ for $\{U_2,T_2\}$ on $H_2$, such that $\psi\oplus\psi_2$ is an orthonormal wavelet for the representation of the group $BS(1,2)$ given by 
$U\oplus U_2$ and $T\oplus T_2$ on $H\oplus H_2$.
\end{theorem}

\begin{proof}
Let 
$$K_2((j,k),(j',k'))=\delta_{j,j'}\delta_{k,k'}-\ip{U^jT^k\psi}{U^{j'}T^{k'}\psi},\quad(j,j',k,k'\in\bz).$$
Then it is easy to check that $K_2$ satisfies \eqref{eqinv1} and \eqref{eqinv2}. Lemma \ref{lemma2} shows that $K_2$ is positive definite. Therefore, by Theorem \ref{thpdtobs}, there exists a Hilbert space $H_2$, a representation $\{U_2,T_2\}$ of the group $BS(1,2)$ and a vector $\psi_2\in H_2$ such that 
$\Span\{U_2^jT_2^k\psi_2\,|\,j,k\in\bz\}=H_2$ and $\ip{U_2^jT_2^k\psi_2}{U_2^{j'}T_2^{k'}\psi_2}=K_2((j,k),(j',k'))$ for all $j,j',k,k'$. 
\par
Then the conclusions follow from Lemma \ref{lemma3}.
\end{proof}

We end this section with a more concrete result which shows that when the Parseval wavelet is {\it semi-orthogonal} we have a more precise description of the orthonormal dilation. 

Recall the definition of a semi-orthogonal wavelet:

\begin{definition}
Let $\psi$ be a Parseval wavelet in some Hilbert space $\mathcal H$ with the representation $\{U,T\}$ of the group $BS(1,2)$. We say that $\psi$ is {\it semi-orthogonal} if $U^jT^k\psi$ is orthogonal to $U^{j'}T^{k'}\psi$ for all $j\neq j'$ in $\bz$ and all $k,k'\in\bz$.
\end{definition}
We will use the notation $\hat U_0=\mathcal FU_0\mathcal F^{-1}$, $\hat T_0=\mathcal FT_0\mathcal F^{-1}$, where 
$\mathcal F$ is the Fourier transform on $\br$, 
$$\mathcal Ff(x)=\int_{\br}f(t)e^{-2\pi ixt}\,dt.$$
\begin{proposition}\label{propsemi}
Any semi-orthogonal Parseval wavelet can be complemented by a Parseval wavelet in a subspace of $L^2(\br)$, and this subspace can be chosen as small as desired. More precisely, let $\psi$ be a semi-orthogonal Parseval wavelet for some representation $U,T$ of the group $BS(1,2)$ on a Hilbert space $H$. Let $\Omega\subset\br$ such that $2\Omega=\Omega$, and $\Omega$ has positive Lebesgue measure. Then there exists a set $F\subset\Omega$ such that, if $\hat\psi_2=\chi_F$, and $\Omega':=\cup_{j\in\bz}2^jF$, then 
$\psi\oplus\psi_2$ is an orthonormal wavelet for $H\oplus\check L^2(\Omega')$ with the representation $\{U\oplus\hat U_0,T\oplus\hat T_0\}$. (Here $\check L^2(\Omega')$ is the Hilbert space of functions in $L^2(\br)$ that have Fourier transform supported on $\Omega'$.)
\end{proposition}

\begin{proof}
Since $\psi$ is a semi-orthogonal Parseval wavelet, $\{T^k\psi\,|\,k\in\bz\}$ is a Parseval frame for its span $W_0$. By \cite{HaLa00}, there exists an isomorphism 
$\mathcal W:W_0\rightarrow L^2[0,1)$ and a subset $E$ of $[0,1)$ such that $\mathcal W\psi=\chi_E$ and $\mathcal WT^kf(x)=e^{2k\pi ix}\mathcal Wf(x)$ for all 
$x\in[0,1)$, $k\in\bz$, and $f\in W_0$. 
\par
By \cite{DaLu96}, there exists a wavelet set $\hat\psi_1=\chi_G$ for $\check L^2(\Omega)$, and this means that the disjoint union $\cup_{j\in\bz}2^jG=\Omega$ and 
$G$ is translation congruent to $[0,1)$, i.e., $\tau: x\mapsto 2x\mod 1$ maps $G$ injectively onto $[0,1)$. Then $G$ has a subset $F$ which is translation congruent to $[0,1)\setminus E$. Of course the sets $2^jF$ will be mutually disjoint for $j\in\bz$. We denote by $\Omega':=\cup_{j\in\bz}2^jF\subset \Omega$.
\par
So $\hat\psi_2=\chi_{F}$ is a Parseval frame for $\check L^2(\Omega')$. Moreover, we have, using the isomorphism $\mathcal W$ and the Fourier transform, for $k\in\bz$:
$$\ip{(T\oplus  T_0)^k(\psi\oplus\psi_2)}{\psi\oplus\psi_2}=\int_0^1e^{2\pi ikx}\chi_E(x)\,dx+\int_{\Omega'}e^{2\pi ikx}\chi_F(x)\,dx=(\ast)$$
and, since $F$ is translation congruent to $[0,1)\setminus E$, (i.e., $\tau$ is injective on $E$)
$$(\ast)=\int_0^1e^{2\pi ikx}\chi_E(x)+\chi_{[0,1)\setminus E}(x)\,dx=\delta_k.$$
The fact that the sets $2^jF$ are mutually disjoint implies that $\hat U_0^j\hat T_0^k\psi_2$ and $\hat U_0^{j'}\hat T_0^{k'}\psi_2$ are orthogonal if $j\neq j'$. Since $\psi$ is semi-orthogonal, the same relation holds for $\psi$, hence it will hold for $\psi\oplus\psi_2$. Consequently, $\{(U\oplus U)^j(T\oplus T)^k\psi\oplus\psi_2\,|\,j,k\in\bz\}$ is an orthonormal family. Since $\psi$ and $\psi_2$ are both Parseval wavelets, it follows using \cite[Proposition 2.5]{HaLa00} that $\psi\oplus\psi_2$ is an orthonormal wavelet for $H\oplus\check L^2(\Omega')$.
\end{proof}

\section{Dilation of MRA Parseval wavelet sets}\label{waveletsets}

We focus now on orthonormal dilations of Parseval wavelet sets. Recall that a Parseval (orthonormal) {\it wavelet set} is a Parseval (orthonormal) wavelet $\psi$ in $L^2(\br)$ such that $\hat\psi=\chi_P$ for some subset $P$ of $\br$. $\hat f$ denotes the Fourier transform of the function $f\in L^1(\br)$:
$$\hat f(x)=\int_{\br}f(t)e^{-2\pi i tx}\,dt,\quad(x\in\br).$$
We restrict our attention to MRA Parseval wavelet sets; we characterize them in Proposition \ref{propm0forwavelet}, and we construct the associated scaling function and low-pass filter. We will see that an orthonormal dilation of a Parseval wavelet set can be realized on a symbolic space, and its precise structure is determined by certain symbolic dynamics (Theorem \ref{thdilmra}). The advantage of this type of orthonormal dilation over the one in Proposition \ref{propsemi} is that the multiresolution structure is preserved too.

We begin with some definitions.

\begin{definition}
The periodization of a function $f$ on $\br$ is
$$\Per(f)(x)=\sum_{k\in\bz}f(x+k),\quad(x\in\br).$$
If $A$ is a subset of $\br$, we denote by 
$\Per(A):=\cup_{k\in\bz}(A+k).$
\end{definition}

\begin{definition}
We will need the following maps
$$\tau(x)=x\mod 1,\quad s(x)=\left(x+\frac12\right)\mod 1, \quad(x\in\br).$$
$$r(x)=x\mod 1,\quad\tau_0(x)=\frac x2,\quad \tau_1(x)=\frac{x+1}{2},\quad(x\in[0,1)).$$
Note that $s(s(x))=x$ for all $x\in[0,1)$.
\end{definition}

\begin{definition}
A subset $A$ of $\br$ is called {\it translation simple} if for all $k\in\bz\setminus\{0\}$, $E\cap (E+k)=\emptyset$ up to Lebesgue measure zero. 
A subset $A$ of $[0,1)$ is called {\it $s$-simple} if it does not contain $x$ and $s(x)$ at the same time, for almost all $x\in [0,1)$, i.e., $A\cap s(A)=\emptyset$ up to measure zero.
\end{definition}

The Parseval wavelet sets are characterized by the following tiling properties:
\begin{proposition}\label{propcharpws} \cite{HaLa00}
Let $\hat\psi=\chi_P$ in $L^2(\br)$. Then $\psi$ is a Parseval wavelet set if and only if $P$ is a multiplicative tile, i.e., $\{2^jP\,|\,j\in\bz\}$ is a partition of $\br$ up to measure zero, and $P$ is translation simple. 
\end{proposition}
\subsection{The scaling function and low-pass filter associated to an MRA Parseval wavelet set}\label{scaling}
 In the next proposition we show how a scaling function and a low-pass filter can be constructed for a MRA wavelet set. For more information on multiresolution analyses see \cite{Dau92}. The wavelet is completely determined by the low-pass filter $m_0=\Per(\chi_M)$: the Fourier transform of the scaling function is an infinite product based on $m_0$ (Proposition \ref{propchardense}), and the wavelet can be obtained from the scaling function by some dilation and translation operations. The orthonormal dilation of the wavelet will be based on the low-pass filter. The function $m_0$ will ``filter'' some symbolic paths, and the dilation will be realized on the set of all the filtered paths.
\begin{proposition}\label{propm0forwavelet}
Let $\psi\in L^2(\br)$ be a Parseval wavelet set, $\hat\psi=\chi_P$. Define 
$$\hat\varphi(x):=\sum_{j\geq 1}\hat\psi(2^jx),\quad(x\in\br^d).$$
Then 
\begin{equation}\label{eqphi}
\hat\varphi=\chi_F, \mbox{ with }F=\cup_{j\geq 1}2^{-j}P,\mbox{ and }F\subset 2F, P=2F\setminus F.
\end{equation}
Assume in addition that $F$ is translation simple. Then there exists a measurable set $M\subset[0,1)$ such that if $m_0=\Per(\chi_M)$, then the following scaling equation is satisfied:
\begin{equation}\label{eqsca}
\hat\varphi(2x)=m_0(x)\hat\varphi(x),\quad(x\in \br),
\end{equation} 
and $m_0$ satisfies the QMF condition 
\begin{equation}\label{eqqmf}
|m_0(x)|^2+|m_0(x+\frac12)|^2=1,\mbox{ i.e., the disjoint union }M\cup s(M)=[0,1).
\end{equation}
Also $m_0$ and $\varphi$ satisfy the conditions in Proposition \ref{propchardense}.
Moreover, in this case, 
\begin{equation}\label{eqwav}
\hat\psi(2x)=(1-m_0(x))\hat\varphi(x),\quad(x\in\br^d),
\end{equation}
i.e., $1-m_0$ is a high-pass filter. 
\end{proposition}
\begin{definition}
If $\hat\psi=\chi_P$ is a Parseval wavelet set such that the set $F$ defined in \eqref{eqphi} is translation simple, we say that $\hat\psi$ is an {\it MRA Parseval wavelet set}.
\end{definition}
\begin{proof}
The relations in \eqref{eqphi} follow directly from the fact that $P$ is a multiplicative tile (Proposition \ref{propcharpws}), therefore the union $F=\cup_{j\geq 1}2^{-j}P$ is disjoint. 
\par
Assume now $F$ is translation simple, so $\tau$ is injective on $F$. Then $\tau(P/2)=\tau(F\setminus F/2)=\tau(F)\setminus\tau(F/2)$. 
Since $F$ is translation simple, $F/2$ is $\frac12\bz$-translation simple, so $\tau(F/2)$ cannot contain both $x$ and $s(x)$ at the same time $(x\in[0,1))$. The same argument works for $\tau(P/2)$. 
\par
Now take $$C:=s\left(\tau(F)\setminus\tau(F/2)\right).$$
Since $\tau(F)\setminus\tau(F/2)=\tau(P/2)$ is translation simple it follows that $C\cap(\tau(F)\setminus\tau(F/2))=\emptyset$ and $C$ is $s$-simple. Moreover $C\cup\tau(F/2)$ is $s$-simple, because $C$ and $\tau(F/2)$ are $s$-simple, and if $x\in C$, then $s(x)\in\tau(F)\setminus\tau(F/2)$ so $s(x)$ is not in $C\cup\tau(F/2)$; if $s(x)\in C$ then $x=s(s(x))$ and we use the same idea.

\par
Since $C\cup\tau(F/2)$ is $s$-simple, we can complete it to an $s$-tile, i.e., there exists a set $D\subset[0,1)$, disjoint from $C\cup\tau(F/2)$ such that if 
$$M:=C\cup\tau(F/2)\cup D$$
then $M$ and $s(M)$ form a partition of $[0,1)$. For example, take 

$$D:=[0,1/2)\setminus\left( \left(\left(C\cup\tau(F/2)\right)\cap[0,1/2)\right)\cup\left(s\left((C\cup\tau(F/2))\cap[1/2,1)\right)\right)\right).$$
Note that $M$ is disjoint from $\tau(F)\setminus\tau(F/2)$. Indeed, $C$ and $\tau(F/2)$ are disjoint from this set, and if $x\in D\cap(\tau(F)\setminus\tau(F/2))$, then
$s(x)\in C$ so $x,s(x)\in M$, a contradiction. 
\par
Since $M$ and $s(M)$ form a partition of $[0,1)$, the function $m_0:=\Per(\chi_M)$ is a QMF filter. We check the scaling equation \eqref{eqsca}: if $m_0(x)\hat\varphi(x)=1$ then $x\in F$ and $x\in\Per(M)$, so $$\tau(x)\in\tau(F)\cap M=(\tau(F/2)\cup(\tau(F)\setminus\tau(F/2)))\,\cap M\subset \tau(F/2).$$
Since $\tau$ is injective on $F$, we get $x\in F/2$ so $\hat\varphi(2x)=1$. 
\par
Conversely, if $\hat\varphi(2x)=1$ then $x\in F/2$ so $\tau(x)\in M$ and $x\in F$, hence $m_0(x)\hat\varphi(x)=1$.

From \eqref{eqsca}, $$(1-m_0(x))\hat\varphi(x)=\hat\varphi(x)-\hat\varphi(2x)=\chi_F(x)-\chi_{F/2}(x)=\chi_{P/2}(x)=\hat\psi(2x).$$
\par
Since $\psi$ is a Parseval wavelet for $L^2(\br)$, $\cup_{j\in\bz}2^jP=\br$ almost everywhere, so condition (ii) in Proposition \ref{propchardense} is satisfied. 
\end{proof}

\subsection{The dilated representation of the Baumslag-Solitar group}\label{dilated}
Proposition \ref{propm0forwavelet} and its proof shows us how to construct the low-pass filter $m_0=\Per(\chi_M)$ associated to our Parseval wavelet set $\hat\psi=\chi_P$. The next step is to construct the representation of the Baumslag-Solitar group $BS(1,2)$ that will contain the dilated orthonormal wavelet. This representation will be supported on a subset of $[0,1)\times\Omega$ where $\Omega$ is the symbolic space
$$\Omega:=\{0,1\}^{\bn}=\{\omega=\omega_1\omega_2\dots\,|\,\omega_n\in\{0,1\}, n\in\bn\}.$$
The subset will be determined by the filter $m_0$. 

\begin{definition}

Let $r:[0,1)\rightarrow[0,1)$, 
$$r(x)=2x\mod 1,\quad(x\in[0,1)).$$
For $x\in[0,1)$, we define $\omega_x\in\{0,1\}$ such that $\tau_{\omega_x}(r(x))=x$. Clearly $\omega_{\tau_kx}=k$ for $k\in\{0,1\}$.
\end{definition}

\begin{definition}
Let $m_0=\chi_M$ be a QMF filter. 
Let $x\in[0,1)$ and $\omega\in\{0,1\}$. We say that the transition $x\rightarrow\tau_{\omega}x$ is {\it possible} if $\tau_\omega x\in M$, i.e., $m_0(\tau_\omega x)=1$. 
Thus if $\tau_\omega x$ is not in $M$, i.e., $m_0(\tau_{\omega}x)=0$, then the transition $x\rightarrow\tau_\omega x$ is not possible.

For each $x\in[0,1)$, because of the QMF equation \eqref{eqqmf}, only one of the transitions $x\rightarrow\tau_0x$, or $x\rightarrow\tau_1x$ is possible. Let $\omega_1\in\{0,1\}$ the digit corresponding to this transition. 
Then for $\tau_{\omega_1}x$ only one of the transitions $\tau_{\omega_1}x\rightarrow\tau_0\tau_{\omega_1}x$ or $\tau_{\omega_1}x\rightarrow\tau_1\tau_{\omega_1}x$ is possible. Let $\omega_2\in\{0,1\}$ be the digit corresponding to this transition. Inductively, there exists a unique $\omega_{n+1}\in\{0,1\}$ such that the transition
$\tau_{\omega_n}\dots \tau_{\omega_1}x\rightarrow\tau_{\omega_{n+1}}\tau_{\omega_n}\dots \tau_{\omega_1}x$ is possible.

We define $\omega(x):=\omega_1\omega_2\dots \in\Omega$ to be the {\it chosen path} for $x$. 

Note that for all $n\geq 1$, $\omega(x)=\omega_1\dots\omega_n\omega(\tau_{\omega_n}\dots\tau_{\omega_1}x)$.

\end{definition}
\begin{remark}
In \cite{CoRa90} a random walk is defined from a low-pass filter $m_0$ on $[0,1)$ with $|m_0(x/2)|^2+|m_0((x+1)/2)|^2=1$. The function $|m_0|^2$ is interpreted as a transition probability. The transition from $x$ to $\tau_ix$ is possible with probability $|m_0(\tau_ix)|^2$ if $m_0(\tau_ix)>0$. We use the same terminology here, however in our case, since $m_0=\Per(\chi_M)$, the walk is actually deterministic.
\end{remark}
\begin{definition}\label{deftilda}
For $x\in[0,1)$ let
$$A(x):=\{\eta\in\Omega\,|\, \eta=\eta_1\dots \eta_n\omega(\tau_{\eta_n}\dots \tau_{\eta_1}x),\mbox{ for some }\eta_1,\dots ,\eta_n\in\{0,1\}\}.$$
Thus, the paths in $A(x)$ start with some random steps $\eta_1,\dots ,\eta_n$, but then follow the chosen path $\omega(\tau_{\eta_n}\dots \tau_{\eta_1}x)$. 
\par
Denote by 
$$\tilde X(m_0):=\{(x,\omega)\,|\,x\in[0,1),\omega\in A(x)\}.$$

Let $$\tilde r:[0,1)\times\Omega\rightarrow[0,1)\times\Omega,\quad \tilde r(x,\omega_1\omega_2\dots)=(r(x),\omega_x\omega_1\omega_2\dots),\quad(x\in[0,1),\omega_1\omega_2\dots\in\Omega).$$
The inverse of $\tilde r$ is
$$\tilde r^{-1}(x,\omega_1\omega_2\dots)=(\tau_{\omega_1}x,\omega_2\omega_3\dots),\quad(x\in[0,1),\omega_1\omega_2\dots\in\Omega).$$
Define the measure $\lambda$ on $[0,1)\times\Omega$ by considering the counting measure on each $A(x)$ and integrating these with respect to $x$ on $[0,1)$:
$$\int_{[0,1)\times\Omega}f\,d\lambda:=\int_0^1\sum_{\omega\in A(x)}f(x,\omega)\,dx.$$
We define the operators $\tilde T$ and $\tilde U$ on $L^2(\tilde X(m_0))$ by:
$$\tilde Tf(x,\omega)=e^{2\pi ix}f(x,\omega),\quad(x\in[0,1),\omega\in\Omega,f\in L^2(\tilde X(m_0))\,),$$
$$\tilde Uf(x,\omega)=\sqrt{2}f(\tilde r(x,\omega)),\quad(x\in[0,1),\omega\in\Omega,f\in L^2(\tilde X(m_0))\,).$$
We define the scaling function 
\begin{equation}\label{eqfitilda}
\tilde\varphi=\chi_{\tilde F}, \mbox{ where }\tilde F=\{(x,\omega(x))\,|\,x\in [0,1)\}.
\end{equation}
Thus the set $\tilde F$ defining the scaling function is obtained by picking exactly the chosen path at each point $x\in[0,1)$.
\end{definition}

\subsection{An encoding of the real numbers}\label{encoding}
We want to realize our dilated representation as a super-representation of the one on $L^2(\br)$. For this we will need to embed $\br$ in the symbolic space $[0,1)\times\Omega$. This will be done by first establishing a one-to-one correspondence between the integers and infinite words that end in either $000\dots$ or $111\dots$. This is the ``two's complement'' encoding system used in computer science, a fact remarked also in \cite{Gun05}. For a more general analysis of this encoding see \cite{DJP07} where it is proved that there are some obstructions when one wants to generalize these encodings to matrix-dilations.
\begin{proposition}\label{corencod}
Let $\underline0$ be  the infinite word $000\dots$ and let $\underline1:=111\dots$.
The map 
$$d_0(\omega_1\dots\omega_n\underline0)=\sum_{k=1}^n\omega_k2^{k-1},(\mbox{ so }d_0(\underline0)=0)$$
is a bijection between $A_0:=\{\omega_1\dots\omega_n\underline 0\,|\,\omega_1,\dots,\omega_n\in\{0,1\}\}$ and $\{k\in\bz\,|\,k\geq 0\}.$
\par
The map 
$$d_1(\omega_1\dots\omega_n\underline 1)=\sum_{k=1}^n\omega_k2^{k-1}-2^{n}, (\mbox{so }d_1(\underline1)=-1)$$
is a bijection between $A_1:=\{\omega_1\dots\omega_n\underline 1\,|\,\omega_1,\dots,\omega_{n}\in\{0,1\}\}$ and $\{k\in\bz\,|\,k<0\}.$
Moreover, for any $\omega\in A_i$, $(i\in\{0,1\})$, and any $x\in[0,1)$,
\begin{equation}\label{eqxpe2n}
\frac{x+d_i(\omega)}{2^n}=\tau_{\omega_n}\dots\tau_{\omega_1}x+d_i(\omega_{n+1}\omega_{n+2}\dots),\quad(n\geq1).
\end{equation}
\end{proposition}

\begin{proof}
The map $d_0$ corresponds to the base $2$ representation of non-negative integers. Note that 
$$d_1(\omega_1\dots\omega_n\underline1)=-(2^n-1-\sum_{k=1}^n\omega_k2^{k-1})-1=-1-\sum_{k=1}^n\breve\omega_k2^{k-1},$$
where $\breve\omega=1-\omega$ for $\omega\in\{0,1\}$. This shows that $d_1$ is also bijective. 
\par
It is enough to prove \eqref{eqxpe2n} for $n=1$, the rest follows by induction. This is obtained by a simple computation (after $n$ steps, one has to use the fact that $d_1(\underline1)=-1$, and $\frac{x-1}{2}=\tau_1x-1$). 
\end{proof}

\begin{definition}
Let $A_\bz:=A_0\cup A_1$. We define the decoding map $d_{\bz}:A_\bz\rightarrow\bz$, 
$$d_\bz(\omega)=\left\{\begin{array}{cc}
d_0(\omega)=\sum_{k=1}^n\omega_k2^{k-1},&\mbox{ if }\omega=\omega_1\dots\omega_n\underline0\in A_0\\
d_1(\omega)=\sum_{k=1}^n\omega_k2^{k-1}-2^{n},&\mbox{ if }\omega=\omega_1\dots\omega_n\underline1\in A_1.\end{array}\right.$$
By Proposition \ref{corencod}, $d_\bz$ is a bijection. 
\par
For each $x\in\br$ define the encoding $\epsilon(x)\in[0,1)\times A_\bz$ as follows: $x$ can be uniquely written as $x=y+k$ with $y\in[0,1)$ and $k\in\bz$, $y:=\tau(x)=x\mod 1, k=x-x\mod1$. Then 
$$\epsilon(x):=(y,d_\bz^{-1}(k))=(x\mod1,d_\bz^{-1}(x-x\mod1)),\quad(x\in\br).$$
\end{definition}

\begin{proposition}
Define the operator $\mathcal E:L^2(\br)\rightarrow L^2([0,1)\times A_\bz,d\lambda)$, 
$$\mathcal E(f)(x,\omega)=f\circ\epsilon^{-1}(x,\omega)=f(x+d_\bz(\omega)).$$ Then $\mathcal E$ is an intertwining isomorphism, $\mathcal E\hat T=\tilde T\mathcal E$, $\mathcal E\hat U=\tilde U\mathcal E$. 
\end{proposition}

\begin{proof}
First we check that $\mathcal E$ is an isometry. This follows from the next computation (for $f\in L^2(\br)$):
$$\int_{\br}|f(x)|^2\,dx=\int_0^1\sum_{k\in\bz}|f(x+k)|^2\,dx=\int_0^1\sum_{\omega\in A_\bz}|f(x+d_\bz(\omega)|^2\,dx.$$
Clearly $\mathcal E$ is invertible, so it is an isomorphism. 
\par
The fact that $\mathcal E$ intertwines the $T$-operators is easy. For the $U$ operators, one only needs to prove that for $x\in[0,1)$ and $\omega\in A_\bz$,
$$2(x+d_\bz(\omega))=r(x)+d_\bz(\omega_x\omega_1\dots),$$
which follows directly from \eqref{eqxpe2n} applied to $r(x)$.
\end{proof}

\subsection{Main result}\label{main}
We can state now the main dilation result of this section:

\begin{theorem}\label{thdilmra}
The operators $\tilde T$ and $\tilde U$ defined in Definition \ref{deftilda} are unitary and $\tilde U\tilde T\tilde U^{-1}=\tilde T^2$. 
A projection $P$ on $L^2(\tilde X(m_0))$ commutes with $\tilde U$ and $\tilde T$ if and only if $P$ is an operator of multiplication by the characteristic function of an $\tilde r$-invariant set, i.e., $Pf=M_{\chi_S}f=\chi_Sf$, where $S\subset [0,1)\times \Omega$ and $\tilde r(S)=S$. 
\par
Let $\tilde\phi=\chi_{\tilde F}$ as in \eqref{eqfitilda}. The translates of $\tilde\varphi$ are orthonormal:
\begin{equation}\label{eqorto}
\ip{\tilde T^k\tilde\varphi}{\tilde\varphi}=\delta_k,\quad(k\in\bz).
\end{equation}
The scaling equation is satisfied:

\begin{equation}\label{eqscal}
\tilde U\tilde\varphi(x,\omega)=\sqrt 2m_0(x)\tilde\varphi(x,\omega),\quad(x\in[0,1),\omega\in\Omega).
\end{equation}
If $\tilde V_0:=\overline{\mbox{span}}\{\tilde T^k\tilde\varphi\,|\,k\in\bz\}$, and $\tilde V_n:=\tilde U^{-n}\tilde V_0$ for $n\in\bz$, then $(\tilde V_n)_{n\in\bz}$ is a multiresolution analysis for $L^2(\tilde X(m_0))$.
\par
Let 
$$\tilde\psi(x,\omega):=\tilde U^{-1}(\sqrt{2}(1-m_0)\tilde\varphi)=\chi_{\tilde P},\quad\mbox{ where }\tilde P=\tilde r(\tilde F)\setminus\tilde F.$$
Then $\tilde\psi$ is an orthonormal wavelet for $L^2(\tilde X(m_0))$. 
\par
Suppose now that $\hat\psi=\chi_P$ is an MRA Parseval wavelet set in $\br$ and let $m_0=\Per(\chi_M)$ be the associated QMF filter, and $\hat\varphi=\chi_F$ be the associated scaling function, as in Proposition \ref{propm0forwavelet}. Then $[0,1)\times A_\bz$ is an $\tilde r$-invariant subset of $\tilde X(m_0)$. Let $\tilde\psi$ be the orthonormal wavelet for $L^2(\tilde X(m_0))$ and let $P_{\br}$ be the corresponding projection $P_\br=M_{\chi_{[0,1)\times A_\bz}}$. Then 
$$P_\br\tilde\varphi=\mathcal E\hat\varphi,\quad P_\br\tilde\psi=\mathcal E\hat\psi.$$

\end{theorem}

\begin{proof}
The operator $\tilde T$ is a multiplication by $e^{2\pi ix}$ so it is unitary. To see that $\tilde U$ is unitary we need the following
\begin{proposition}
For all integrable functions $f$ on $[0,1)\times \Omega$,
\begin{equation}\label{eqinme}
\int_{[0,1)\times\Omega}2f\circ\tilde r\,d\lambda=\int_{[0,1)\times\Omega}f\,d\lambda.
\end{equation}
\end{proposition}

\begin{proof}
The Lebesgue measure on $[0,1)$ has the following strong invariance property:

\begin{equation}\label{eqstrinv}
\int_0^1f(x)\,dx=\int_0^1\frac{1}{2}\sum_{\omega\in\{0,1\}}f(\tau_\omega x)\,dx,\quad(f\in L^1[0,1)).
\end{equation}
Using equation \eqref{eqstrinv} we have:
$$\int_{[0,1)\times\Omega}2f\circ\tilde r\,d\lambda=\int_0^12\sum_{\omega\in A(x)}f(r(x),\omega_x\omega)\,dx=\int_0^1\sum_{k\in\{0,1\}}\sum_{\omega\in  A(\tau_kx)}f(r(\tau_kx),\omega_{\tau_kx}\omega)\,dx=$$
$$\int_0^1\sum_{k\in\{0,1\}}\sum_{\omega\in A(\tau_kx)}f(x,k\omega)\,dx=\int_0^1\sum_{\omega\in A(x)}f(x,\omega)\,dx,$$
and, for the last equality, we used the fact that $A(x)$ is the disjoint union 
\begin{equation}\label{eqax}
A(x)=0A(\tau_0x)\cup 1A(\tau_1x).
\end{equation}
This proves \eqref{eqinme}.

\end{proof}

 Equation \eqref{eqinme} shows that $\tilde U$ is an isometry and since $\tilde r$ is bijective and the set $\tilde X(m_0)$ is invariant under $\tilde r$, the operator $U$ is unitary. 
\par
The relation $\tilde U\tilde T\tilde U^{-1}=\tilde T^2$ is obtained by an easy computation. 
\par
Let $W$ be a projection that commutes with $\tilde U$ and $\tilde T$. Then $W$ commutes with $\sum_ka_k\tilde T^k$. So it must commute with all operators of multiplication by functions that depend only on $x$, $M_gf(x,\omega)=g(x)f(x,\omega)$. Then $W$ commutes with operators of the form $\tilde U^{-n}M_g\tilde U^{n}$, $n\in\bn$, but these are operators of multiplication by $g\circ\tilde r^{-n}$, i.e., operators of multiplication by functions which depend only on $x$ and $\omega_1,\dots,\omega_n$. The SOT-closure of these operators is the algebra of all multiplication operators on $L^2(\tilde X(m_0))$. But this is a maximal abelian algebra, so $W$ must be contained in it. Thus $W$ is a multiplication operator $W=M_f$. Since $W$ is a projection, $f$ is a characteristic function $f=\chi_S$. Since $W$ commutes with $\tilde U$, the set $S$ is $\tilde r$-invariant. 
\par
The orthogonality of the translates of $\tilde\varphi$ is trivial. The scaling equation follows from the following equality:
$$\tilde r^{-1}(\tilde F)=\{(x,\omega)\,|\, \omega(r(x))=\omega_x\omega_1\dots\}=\{(x,\omega)\,|\,\tau_{\omega_x}(r(x))\in M,\omega(x)=\omega_1\dots\}=(M\times\Omega)\cap\tilde F.$$
\par
Notice that $\tilde V_0$ consists of all the functions supported on $\tilde F$. Then $\tilde V_n$ consists of the functions supported on $\tilde r^{n}\tilde F$ for all $n\in\bz$. We also have $\tilde r^{-1}\tilde F\subset\tilde F$, from the scaling equation. The definition of $A(x)$ implies that $\cup_{n\geq0}\tilde r^n\tilde F=\tilde X(m_0)$. This implies that the union of the subspaces $\tilde V_n$ is dense. 
\par
To show that $\cap_{n}V_n=\{0\}$, note that, using \eqref{eqinme}
$$\lambda(\tilde r^{-n}\tilde F)=\int_{[0,1)\times\Omega}\chi_{\tilde F}\circ\tilde r^n\,d\lambda=\frac{1}{2^n}\int_{[0,1)\times\Omega}\chi_{\tilde F}\,d\lambda=\frac{1}{2^n}.$$
Therefore a function $f$ in $\cap V_n$ is supported on a set of measure $0$, so it has to be identically $0$.

Now let us consider the case of a MRA Parseval wavelet set. Proposition \ref{propchardense} shows that $A(x)$ contains $A_\bz$ for almost every $x\in[0,1)$. We have that $P_{\br}\tilde\varphi=M_{\chi_{[0,1)\times A_\bz}}\chi_{\tilde F}=\chi_{[0,1)\times A_\bz\,\cap\tilde F}$. We have to check that 
\begin{equation}\label{eqcap}
([0,1)\times A_\bz)\cap\tilde F=\epsilon(F).
\end{equation}
If $x\in[0,1)$ and $\omega\in A(x)\cap A_\bz$ then the chosen path of $x$ ends in $\underline0$, in which case we let $i:=0$, or $\underline1$, and in this case we let $i:=1$. Since $\omega$ is the chosen path of $x$, using equation \eqref{eqxpe2n}, we have that $m_0((x+d_i(\omega))/2^n)=m_0(\tau_{\omega_n}\dots\tau_{\omega_1}x)=1$. Therefore, since $\hat\varphi$ is the infinite product in Proposition \ref{propchardense}, we obtain that $\hat\varphi(x+d_i(\omega))=1$ so $(x,\omega)\in \epsilon(F)$. 

Conversely, if $(x,\omega)\in \epsilon(F)$, then $x+d_i(\omega)$ is in $F$, where $i$ is $0$ if $\omega$ ends in $\underline0$ and $i=1$ if $\omega$ ends in $\underline1$. So $m_0((x+d_i)/2^n)=1$ and with equation \eqref{eqxpe2n}, this shows that $m_0(\tau_{\omega_n}\dots\tau_{\omega_1}x)=1$, which implies that $\omega$ is the chosen path of $x$. This proves \eqref{eqcap}. Since $\mathcal E$ and $P_{\br}$ intertwine the representations, and since the relation between $\tilde\psi$ and $\tilde\varphi$ is the same as the relation between $\hat\psi$ and $\hat\varphi$, it follows that $P_{\br}\tilde\psi=\mathcal E\hat\psi$.
\end{proof}

The next proposition characterizes the density property of the multiresolution in terms of the chosen paths. It will tell us in which cases the orthonormal dilation contains $L^2(\br)$ as a subrepresentation.

\begin{proposition}\label{propchardense}
Suppose $\hat\varphi=\chi_F$ and $m_0=\Per(\chi_{M})$ where $F\subset\br$, $M\subset [0,1)$. Assume that:
$$\hat\varphi(x)=\prod_{n=1}^\infty m_0\left(\frac{x}{2^n}\right),\mbox{ i.e., }F=\bigcap_{n=1}^\infty2^n\Per(M).$$
The following affirmations are equivalent:
\begin{enumerate}
\item $\lim_{n\rightarrow\infty}\hat\varphi\left(\frac{x}{2^n}\right)=1$ for a.e. $x\in\br$;
\item $\cup_{n=1}^\infty 2^nF=\br$ (up to measure zero);
\item $\lim_{n\rightarrow\infty}m_0\left(\frac{x}{2^n}\right)=1$ for a.e., $x\in\br$;
\item $\lim_{n\rightarrow\infty}m_0(\tau_0^n{x})=1$ and  $\lim_{n\rightarrow\infty}m_0(\tau_1^n{x})=1$, for a.e. $x\in[0,1)$;
\item For a.e. $x\in[0,1)$, the chosen paths  $\omega(\tau_0^n{x})=\underline 0$ and $\omega(\tau_1^n{x})=\underline 1$ if $n$ is big enough;
\item For a.e. $x\in[0,1)$, the set $A(x)$ contains the paths $\omega_1\dots\omega_n\underline 0$ and $\omega_1\dots\omega_n\underline 1$ for all $\omega_1,\dots,\omega_n\in\{0,1\}$.
\end{enumerate}
\end{proposition}

\begin{proof}
(i)$\Rightarrow$(ii). If $\hat\varphi(x/2^n)\rightarrow 1$, then $x/2^n\in F$ for $n$ big enough, so $x\in 2^nF$. 
\par
(ii)$\Rightarrow$(iii). If $x\in 2^nF$ then, since $F\subset 2F$ (from the hypothesis), $x/2^{n+k}\in F$ for $k\geq 0$. But $F\subset 2\Per(M)$, so
$m_0(x/2^{n+k})=1$ for $k\geq 1$. 
\par
(iii)$\Rightarrow$(iv). For any $x\in[0,1)$ and any $k\in\bz$ we have $m_0((x+k)/2^n)=1$ for $n$ big enough. Using the encoding in Corollary \ref{corencod}, we obtain that $m_0(\tau_0^k\tau_{\omega_n}\dots\tau_{\omega_1}{x})=1$ and $m_0(\tau_1^k\tau_{\omega_n}\dots\tau_{\omega_1}{x})=1$ if $k$ is big enough, for all $\omega_1,\dots\omega_n\in\{0,1\}$. (iv) is a particular case of this. 
\par
(iv)$\Rightarrow$(v). Evident. 
\par
(v)$\Rightarrow$(vi). Let $\omega_1,\dots\omega_n\in\{0,1\}$. Then apply (v) to $\tau_{\omega_n}\dots\tau_{\omega_1}x$ and (vi) follows. 
\par
(vi)$\Rightarrow$(i). We have $m_0(\tau_0^k\tau_{\omega_n}\dots\tau_{\omega_1}x)=1$ for all $\omega_1,\dots,\omega_n$ and $k$ big enough. Similarly with $\tau_1$. Using the encoding in Corollary \ref{corencod}, we obtain $m_0((x+k)/2^n)=1$ for all $x\in[0,1)$ and all $k\in\bz$, and $n$ big enough. But this implies that for all $x\in\br$, $x/2^n\in\Per(M)$ for $n$ big enough, so $x/2^p\in F$ for some $p\geq1$. 
\end{proof}

\subsection{Cyclic paths}\label{cyclic}
Our construction of the orthonormal dilation is based on finding the chosen paths. We will show that under some extra assumption on $m_0$ the chosen paths are eventually periodic, and the orthonormal dilation has a particularly simple form and can be realized on an orthogonal sum of copies of $L^2(\br)$ just as in \cite{BDP05}.
\begin{definition}
We call a set $C:=\{\theta_0,\dots,\theta_{p-1}\}$ in $[0,1)$ a cycle corresponding to $l_0\dots l_{p-1}\in\{0,1\}^p$, if $\tau_{l_0}\theta_0=\theta_1, \tau_{l_1}\theta_1=\theta_2,\dots, \tau_{l_{p-2}}\theta_{p-2}=\theta_{p-1}$ and $\tau_{l_{p-1}}\theta_{p-1}=\theta_0$. 

We denote by $\underline{l_0\dots l_{p-1}}$ the infinite word obtained by the infinite repetition of the finite word $l_0\dots l_{p-1}$, i.e., 
$$\underline{l_0\dots l_{p-1}}=l_0\dots l_{p-1}l_0\dots l_{p-1}\dots .$$
We denote by $\Omega_C$ the set of infinite words that end in $\underline{l_0\dots l_{p-1}}$, i.e.,
$$\Omega_C:=\{\omega_1\dots\omega_n\underline{l_0\dots l_{p-1}}\,|\,\omega_1,\dots,\omega_n\in\{0,1\}\}.$$
\end{definition}

We define the encoding/decoding maps between eventually cyclic paths and integers as in \cite{DJP07}.
\begin{definition}
Let $C=\{\theta_0,\dots,\theta_{p-1}\}$ be a cycle, $\bz_p:=\{0,1,\dots,p-1\}$. Let $T_0$ and $U_0$ be the operators on $L^2(\br)$ from \eqref{eqwv2}. Define the following operators on $L^2(\br\times\bz_p)$:
\begin{equation}\label{eqtc}
\hat T_C(f_0,\dots,f_{p-1})=(e^{2\pi i\theta_0}\hat T_0f_0,\dots,e^{2\pi i\theta_{p-1}}\hat T_0f_{p-1}),
\end{equation}

\begin{equation}\label{equc}
\hat U_C(f_0,\dots,f_{p-1})=(\hat U_0f_{p-1},\hat U_0f_0,\dots,\hat U_0f_{p-2}).
\end{equation}
Note that equation \eqref{equc} can be rewritten as
\begin{equation}\label{equc2}
\hat U_Cf=\sqrt2f\circ\alpha_C,\quad f\in L^2(\br\times\bz_p),\mbox{ where }\alpha_C(x,j)=(2x,(j-1)\mod p),\quad(x\in\br, j\in\bz_p).
\end{equation}
We define the decoding map
\begin{equation}\label{eqe_C}
d_C:[0,1)\times\Omega_C\rightarrow\br\times\bz_p,\quad d_C(x,\omega)=(x-\theta_{j(\omega)}+k(\omega),j(\omega)),
\end{equation}
where $j(\omega)\in\bz_p$ and $k(\omega)\in\bz$ are defined as follows: there is a unique $j(\omega)\in\{0,\dots,p-1\}$ such that $\omega=\omega_0\dots\omega_{np-1}\underline{l_{j(\omega)}l_{j(\omega)+1}\dots l_{p-1}l_0\dots l_{j(\omega)-1}}$ for some $\omega_1,\dots,\omega_{np-1}\in\{0,1\}$. 
\begin{equation}\label{eqkomega}
k(\omega)=\omega_0+\dots+2^{np-1}\omega_{np-1}+\theta_{j(\omega)}-2^{np}\theta_{j(\omega)}.
\end{equation}

\end{definition}

\begin{remark}
The inverse transformation $\br\times\bz_p\ni(x,j)\mapsto (y,\omega)\in[0,1)\times\Omega_C$ is constructed as follows:
There is a unique $y\in[0,1)$ and a $k\in\bz$ such that $x-\theta_j=y+k$.

We will associate to $(k,j)$ a path $\omega\in\Omega_C$. The way to define $\omega$ resembles the Euclidian algorithm. 

First we define the map $\mathcal R_C:\bz-C\rightarrow\bz-C$, using a division with remainder: for $a-\theta_j\in\bz-\theta_j$ there is a unique $\mathcal R_C(a-\theta_j)\in\bz-\theta_{j+1}$ and $d\in \{0,1\}$ such that 
$$a-\theta_{j}=2\mathcal R_C(a-\theta_j)+d.$$
(Note that we use here the notation $\theta_j=\theta_{j\mod p}$.)

Then, to define $\omega\in\Omega_C$ from $k\in\bz$ and $j\in\bz_p$, we iterate this division and keep the remainders:
there is a unique $\omega_0\in\{0,1\}$ such that $k-\theta_j=2\mathcal R_C(k-\theta_j)+\omega_0$; at the next step, there is a unique $\omega_1\in\{0,1\}$ such that $\mathcal R_C(k-\theta_j)=2\mathcal R_C^2(k-\theta_j)+\omega_1$; at step $n$, there is a unique $\omega_n\in\{0,1\}$ such that $$\mathcal R_C^n(k-\theta_j)=2\mathcal R_C^{n+1}(k-\theta_j)+\omega_n.$$
Then $\omega$ is defined by $\omega_0\omega_1\dots$.
\end{remark}

Using the decoding maps, one can embed the represenatation associated to a cycle into the representation of the group $BS(1,2)$ on the symbolic space $L^2([0,1)\times\Omega)$ defined in Section \ref{dilated}.

\begin{theorem}\label{thencc}\cite{DJP07}
The map $d_C$ is bijective and 
\begin{equation}\label{eqecr}
d_C\circ\tilde r=\alpha_C\circ d_C
\end{equation}
The map $\mathcal E_C: L^2(\br\times\bz_p)\rightarrow L^2([0,1)\times\Omega_C,\lambda)$, $\mathcal E_Cf=f\circ d_C$ is an isometric isomorphism that intertwines the representations $\{\hat U_C,\hat T_C\}$ and $\{\tilde U, \tilde T\}$.
\end{theorem}

\begin{proposition}\label{propc} Let $C$ be the cycle corresponding to $l_0\dots l_{p-1}$. Let
$$m_0^{(p)}(x)=m_0(x)m_0(rx)\dots m_0(r^{p-1}x)=m_0(x)m_0(2x)\dots m_0(2^{p-1}x),\quad(x\in\br).$$
 The following affirmations are equivalent:
\begin{enumerate}
\item For a.e. $x\in[0,1)$, $\lim_{n\rightarrow\infty} m_0^{(p)}((\tau_{l_{p-1}}\dots\tau_{l_0})^nx)=1$;
\item For a.e. $x\in[0,1)$, $A(x)\supset\Omega_C$;
\item The representation $\pi_{m_0}:=\{\tilde U,\tilde T\}$ on $L^2(\tilde X(m_0))$ contains $\pi_C=\{\hat U_C,\hat T_C\}$ as a subrepresentation. 
\end{enumerate}
\end{proposition}

\begin{proof}
(i)$\Rightarrow$(ii) Since the set of finite words is countable, and since the maps $\tau_\omega$ and $x\mapsto 2x\mod 1$ preserve sets of measure zero, we have that for a.e. $x\in[0,1)$, $\lim_{n\rightarrow\infty}m_0^{(p)}((\tau_{l_{p-1}}\dots\tau_{l_0})^n(\tau_{\omega_m}\dots\tau_{\omega_1}x))=1$ for all $\omega_1,\dots,\omega_m\in\{0,1\}$. But this means that, if $\omega=\omega_1\dots\omega_m\underline{l_0\dots l_{p-1}}$ then $m_0(\tau_{\omega_n}\dots\tau_{\omega_1}x)=1$ for $n$ large enough. And this implies that if we choose $n$ large, the chosen path of  $\tau_{\omega_{n}}\dots\tau_{\omega_1}x$ is $\omega_{n+1}\omega_{n+2}\dots$. Thus any such $\omega$ is in $A(x)$ which implies (ii)

(ii)$\Rightarrow$(iii) is clear from Theorem \ref{thencc}.

(iii)$\Rightarrow$(ii) We need a lemma:

\begin{lemma}\label{lem}
Let $A$ be a map from $[0,1)$ to countable subsets of $\Omega$. Assume that 
$$\tilde X(A):=\{(x,\omega)\,|\,x\in[0,1),\omega\in A(x)\},$$
is invariant under $\tilde r$. 
Consider representations of the form $\pi_A:=\{\tilde U,\tilde T\}$ on $L^2(\tilde X(A),\lambda)$. If $A_1$ and $A_2$ are such maps, then $\pi_{A_1}$ is a subrepresentation of $\pi_{A_2}$ if and only if $A_1(x)\subset A_2(x)$ for almost every $x\in[0,1)$.
\end{lemma}

\begin{proof}
The sufficiency is immediate. 
Let $W$ be an isometry between $L^2(\tilde X(A_1))$ to $L^2(\tilde X(A_2))$ that intertwines the representations. Then, proceeding as in the proof of Theorem \ref{thdilmra}, $W$ must intertwine multiplication operators on the two spaces $\tilde X(A_1)$ and $\tilde X(A_2)$. Therefore $\tilde X(A_1)\cap\tilde X(A_2)$ cannot be empty and since $W$ is an isometry we must have $\tilde X(A_1)\subset \tilde X(A_2)$. The lemma follows.
\end{proof}

Since $\pi_{m_0}$ contains $\pi_C$, using Lemma \ref{lem}, we have that $\tilde X(m_0)$ must contain $[0,1)\times\Omega_C$, and this implies (ii).

(ii)$\Rightarrow$(i) We have that for a.e. $x\in[0,1)$, $\omega_1\omega_2\dots:=\underline{l_0\dots l_{p-1}}$ is in $A(x)$. So for some $n$, the chosen path of $\tau_{\omega_n}\dots\tau_{\omega_1}x$ is $\omega_{n+1}\omega_{n+2}\dots$. We make $n$ bigger if necessary to have $n$ of the form $n=kp$. But his implies that $m_0(\tau_{l_0}(\tau_{l_{p-1}}\dots\tau_{l_0})^kx)=1$, $m_0(\tau_{l_1}\tau_{l_0}(\tau_{l_{p-1}}\dots\tau_{l_0})^kx)=1$, and so on. Therefore $m_0^{(p)}((\tau_{l_{p-1}}\dots\tau_{l_0})^mx)=1$ for $m$ large enough. And this implies (i).
\end{proof}

\begin{definition} Let $m_0=\Per(\chi_M)$ be a QMF filter.
\begin{enumerate}
 \item We call $M$ (and $m_0$) {\it partitionable} if there exists a finite partition $I_1,\dots, I_q$ of $M$ with the property that for each $i\in\{1,\dots, q\}$ there exists a $j(i)\in\{1,\dots,q\}$ and a $\nu(i)\in\{0,1\}$ such that $\tau_{\nu(i)}(I_i)\subset I_{j(i)}$. We say that the partition $(I_i)_{i=1}^q$ is {\it subordinated} to $M$ (and $m_0$). 
\item
For the partition $(I_i)_{i=1}^q$, we construct the following graph: the vertices are the intervals $I_i$, $i\in\{1,\dots,q\}$. We have an edge from $i$ to $j$ if and only if $j=j(i)$; moreover we label the edge from $i$ to $j(i)$ by $\nu(i)$. We call this {\it the graph associated to the partition} $(I_i)_{i=1}^q$
\item
For each cycle in the graph associated to the partition $(I_i)_{i=1}^q$, let $l_0\dots l_{p-1}$ be the corresponding labels. We say that the cycle $C$ associated to the word $l_0\dots l_{p-1}$ is a cycle {\it associated to the partition} $(I_i)_{i=1}^q$.
\end{enumerate}
\end{definition}

\begin{theorem}\label{th3_24}
Let $m_0=\Per(\chi_M)$ be a partitionable QMF filter, and let $(I_i)_{i=1}^q$ be a partition subordinated to $m_0$. 
\begin{enumerate}
\item The representation $\pi_{m_0}=\{\tilde U,\tilde T\}$ on $L^2(\tilde X(m_0))$ is a subrepresentation of $$\oplus\{\pi_C\,|\, C\mbox{ cycle associated to the partition }(I_i)_{i=1}^q\}.$$
(Recall $\pi_C=\{\hat U_C,\hat T_C\}$ on $\underbrace{L^2(\br)\oplus\dots L^2(\br)}_{ \mbox{length}(C)\mbox{-times}}$.)
\item If in addition all cycles $C$ associated to the partition $(I_i)_{i=1}^q$ are contained in the interior of $M$, then 
$$\pi_{m_0}=\oplus\{\pi_C\,|\, C\mbox{ cycle associated to the partition }(I_i)_{i=1}^q\}.$$
\end{enumerate}
\end{theorem}

\begin{proof}(i)
We will show that for a.e. $x\in[0,1)$, $A(x)\subset\cup\Omega_C$ where the union is done over all the cycles associated to the partition. 

Take $x\in[0,1)$, and let $\omega=\omega_1\omega_2\dots$ be its chosen path. Then $\tau_{\omega_1}x\in M$, so there is some $i_0\in\{1,\dots,q\}$ such that $\tau_{\omega_1}x\in I_{i_0}$. Also $\tau_{\omega_2}\tau_{\omega_1}x\in M$, but, since $\tau_{\omega_1}x\in I_{i_0}$, this implies that $\tau_{\omega_2}\tau_{\omega_1}x\in I_{j(i_0)}$ and $\omega_2=\nu(i_0)$. By induction, we obtain $\omega_3=\nu(j(i_0)), \dots, \omega_{n}=\nu(j^{n-2}(i_0))$, where $j^n=j\circ\dots\circ j$, $n$ times. Moreover, we have that $\tau_{\omega_{n+1}}\dots\tau_{\omega_1}x\in I_{j^n(i_0)}$, so $\nu^n(i_0)$ is the label for the edge between $j^{n-1}(i_0)$ and $j^n(i_0)$. Since the graph is finite it is clear that this procedure will enter a cycle, i.e., the sequence $j^n(i_0)$ and $\nu^n(i_0)$ are eventually periodic. The cycle is associated to the partition, and this proves that the chosen path $\omega$ of $x$ is in one of the sets $\Omega_C$. From this it follows immediately that $A(x)$ is contained in the union of the the sets $\Omega_C$. Then (i) follows, since $\tilde X(m_0)$ is subset of $\cup\Omega_C$ which is invariant under $\tilde r$.

(ii) We use Proposition \ref{propc}. We have that all cycles associated to the partition are interior points for $M$. Let $C=\{\theta_0,\dots,\theta_{p-1}\}$ be such a cycle and let $l_0\dots l_{p-1}$ be the corresponding word. We have that $(\tau_{l_{p-1}}\dots\tau_{l_0})^nx$ converges to the fixed point of the map $\tau_{l_{p-1}}\dots\tau_{l_0}$ which is $\theta_0$. Therefore $m_0((\tau_{l_{p-1}}\dots\tau_{l_0})^nx)=1$ for $n$ large enough. Similarly for the other cyclic permutations $l_1\dots l_{p-1}l_0$ and so on. This implies that $m_0^{(p)}((\tau_{l_{p-1}}\dots\tau_{l_0})^nx)=1$ for $n$ large enough, and with Proposition \ref{propc}, we get (ii).
\end{proof}

\section{Examples}\label{ex}
\begin{example}\label{ex1}
Consider 
$$\hat\psi:=\chi_{[-2a,-a]\cup[a,2a]},\quad (0<a\leq\frac14).$$
Since $P:=[-2a,a]\cup[a,2a]$ is a dilation tile, and translation simple, $\psi$ is a Parseval wavelet. We want to use our theory to construct an orthonormal dilation. We will see that:

\begin{proposition}
 The wavelets $\hat\psi=\chi_{[-2a,-a]\cup[a,2a]}$, $0<a\leq\frac14$ have an orthonormal dilation in the space 
$L^2(\br)\oplus L^2(\br)\oplus L^2(\br)$ with the representation $\{U_0\oplus \hat U_C,T_0\oplus \hat T_C\}$, where $C$ is the cycle $C:=\{\frac13,\frac23\}$.
\end{proposition}

First, we have to compute the associated scaling function and low-pass filter. By Proposition \ref{propm0forwavelet}, we have that the scaling function is $\hat\varphi=\chi_F$, with $F=\cup_{j\geq1}2^{-j}P=[-a,a]$. This set is translation simple.

To construct the set $M$ for the low-pass filter, we follow the procedure in the proof of Proposition \ref{propm0forwavelet}. Recall $\tau (x)=x\mod1$, $s(x)=(x+1/2)\mod 1$. We have
$$\tau(F)=[0,a]\cup[1-a,1],\quad\tau(F/2)=[0,\frac a2]\cup[1-\frac a2,1].$$
Then 
$$C:=s(\tau(F)\setminus\tau(F/2))=s([\frac a2,a]\cup [1-a,1-\frac a2])=[\frac12+\frac a2,\frac12+a]\cup[\frac12-a,\frac12-\frac a2].$$
The set $M$ must contain both sets $\tau(F/2)$ and $C$, and it must be disjoint from the sets $s(C)$ and $s(\tau(F/2))$.

We are left with an ``undecided zone'', $[a,1/2-a]\cup[1/2+a,1-a]$, where we must make a choice of a subset $D$ with the property that $|\{x,s(x)\}\cap D|=1$ for all $x$ in this zone. Note that $s$ maps the two intervals of this zone into each other. 

We {\it pick} here 
$$D:=[\frac14,\frac 12-a]\cup[\frac 12+a,\frac34].$$
Of course there are many other choices, and it would be interesting to see how these choices will affect the dilation.

Then we get that the support set for our low-pass filter is 
$$M:=[0,\frac a2]\cup[\frac14,\frac12-\frac a2]\cup[\frac12+\frac a2,\frac34]\cup[1-\frac a2,1].$$

Next we have to see what the chosen paths are. For this we find a partition subordinated to $M$. This is easy. The four intervals will give us this partition. Indeed we have that 
$$\tau_0[0,\frac a2]\subset [0,\frac a2],\quad \tau_1[1-\frac a2,1]\subset [1-\frac a2,1],$$
$$\tau_1[\frac14,\frac12-\frac a2]\subset [\frac 12+\frac a2,\frac34],\quad \tau_0[\frac12+\frac a2,\frac34]\subset[\frac14,\frac12-\frac a2].$$

Therefore we have the following cycles associated to the partition: $\underline 0$, $\underline 1$ (the occurence of these two cycles should be no surprise because our filter comes from a construction in $\br$, where the low-pass condition on $\chi_M$ implies that these cycles are present), and $\underline{10}$ (or $\underline{01}$). The cycle corresponding to $\underline 0$ is $0$, the one for $\underline1$ is $1$, the cycle for $\underline{10}$ is $c:=\{\frac13,\frac23\}$. So 
$\theta_0=\frac13, \theta_1=\frac23$, and since $\tau_1\frac13=\frac23, \tau_0\frac23=\frac13$, we have $l_0=1$ and $l_1=0$. Note also that the cycle $\{\frac13,\frac23\}$ lies in the interior of $M$.

Since we have these cycles, our dilation will be constructed in the space $[0,1)\times \{\omega\in\Omega\,|\, \omega \mbox{ ends in }\underline0,\underline1\mbox{ or }\underline{10}\}$. Or equivalently, using the encoding/decoding, it can be done in $\br\times\{*,0,1\}$, where $*$ will be the index for the $L^2(\br)$-component that we started from (corresponding to $\underline 0,\underline1$), and the other two components $\{0,1\}$ will correspond to the cycle $\underline{10}$.

Next, we want to find what the dilated scaling function $\tilde\varphi=\chi_{\tilde F}$ is, so we have to find the set $\tilde F$. Recall that 
$$\tilde F=\{(x,\omega(x))\,|\,x\in[0,1)\},$$
where $\omega(x)$ is the chosen path of $x$. 

To determine the chosen path for a point $x\in[0,1)$ we actually need to find only the first digit, i.e., to find $\omega_1\in\{0,1\}$ such that $\tau_{\omega_1}x\in M$, because once $\tau_{\omega_1}x$ is in $M$, we use the partition associated to $M$ to see what the next digits of the chosen path are. Using this, and the fact that $0<a<1/4$ we obtain:
$$\omega(x)=\left\{\begin{array}{cc}
\underline0,&\mbox{ if }x\in[0,a]\\
\underline{10},&\mbox{ if }x\in[a,\frac12]\\
\underline{01},&\mbox{ if }x\in[\frac12,1-a]\\
\underline{1},&\mbox{ if }x\in[1-a,1].
\end{array}\right.$$

Now that we have the chosen paths for each $x$ in $[0,1)$ we use the decoding maps to see how the set $\tilde F$ is mapped inside $\br\times\{*,0,1\}$.

On $[0,a]$ we have $\omega(x)=\underline 0$. Then $d_0(\underline 0)=0$ so $\epsilon^{-1}([0,a]\times\{\underline0\})=[0,a]+0=[0,a]$.

On $[1-a,1]$ we have $\omega(x)=\underline1$. Then $d_1(\underline 1)=-1$ so 
$\epsilon^{-1}([1-a,1]\times\{\underline1\})=[1-a,1]-1=[-a,0]$.

Therefore the first component of the set (the one corresponding to $*$) will be $[-a,0]\cup[0,a]=[-a,a]$. This is to be expected, of course, since that was our objective: to dilate the wavelet and scaling function that we started with, so the $*$ component of the scaling function $\chi_{\tilde F}$ should be $\chi_F$. Similarly for the wavelet $\tilde\psi$.
On $[a,\frac12]$ we have $\omega(x)=\underline{10}$. Then $\theta_{j(\underline{10})}$ is the fixed point of $\tau_0\tau_1$, which is $\frac13=\theta_0$ so $j(\underline{10})=0$. Then $k(\underline{10})=0$, and $d_c(x,\underline{10})=(x-\frac13,0)$. So $d_c([a,\frac12]\times\{\underline{10}\})=[a-\frac13,\frac16]\times\{0\}$.

On $[\frac12,1-a]$ we have $\omega(x)=\underline{01}$. Then $\theta_{j(\underline{01})}$ is the fixed point of $\tau_1\tau_0$, which is $\frac23=\theta_1$ so $j(\underline{01})=1$. Then $k(\underline{01})=0$, and $d_c(x,\underline{01})=(x-\frac23,0)$. So $d_c([\frac12,1-a]\times\{\underline{01}\})=[-\frac16,\frac13-a]\times\{1\}$.

Consequently we have that 
$$\tilde\varphi\circ (\epsilon,d_c^{-1})=(\chi_{[-a,a]},\chi_{[a-\frac13,\frac16]},\chi_{[-\frac16,\frac13-a]}).$$

Let $\alpha(x_*,x_0,x_1)=(2x_*\mod1, 2x_1\mod1,2x_0\mod_1)$. The support set of the dilated wavelet is $\tilde P=\alpha(\tilde F)\setminus \tilde F$.
We have $\alpha(\tilde F)=[-2a,2a]\times\{*\}\cup[-1/3,2/3-2a]\times\{0\}\cup[2a-2/3,1/3]\times\{1\}$. Therefore
$$\tilde P=([-2a,-a]\cup[a,2a])\times\{*\}\cup([-\frac13,a-\frac13]\cup[\frac16,\frac23-2a])\times\{0\}\cup([2a-\frac23,-\frac16]\cup[\frac13-a,\frac13])\times\{1\}.$$
Finally the dilated orthonormal wavelet is
$$\tilde\psi=(\chi_{[-2a,-a]\cup[a,2a]},\chi_{[-\frac13,a-\frac13]\cup[\frac16,\frac23-2a]},\chi_{[2a-\frac23,-\frac16]\cup[\frac13-a,\frac13]}).$$
\end{example}

\begin{example}\label{ex2}
Let us consider again the wavelet set in the previous example, now with $a=\frac18$. So $\hat\psi=\chi_{[-\frac14,-\frac18]\cup[\frac18,\frac14]}$.
\begin{proposition}
The wavelet $\hat\psi=\chi_{[-\frac14,-\frac18]\cup[\frac18,\frac14]}$ has an orthonormal dilation in the space $L^2(\br)\oplus L^2(\br)\oplus L^2(\br)\oplus L^2(\br)$ with the representation $\{U_0\oplus \hat U_C,T_0\oplus \hat T_C\}$ where $C$ is the cycle $C:=\{\frac17,\frac47,\frac27\}$. This proves that the representation associated to orthonormal dilations is not unique.
\end{proposition}
 We saw that for this wavelet set we can take the scaling function to be $\hat\varphi=\chi_{[-\frac18,\frac18]}$. The support $M$ for the low-pass filter must contain $[0,\frac1{16}]\cup[\frac38,\frac7{16}]\cup[\frac9{16},\frac58]\cup[\frac{15}{16},1]$, and it should be disjoint
from $[\frac1{16},\frac18]\cup[\frac7{16},\frac9{16}]\cup[\frac78,\frac{15}{16}]$. And we have the undecided zone
$[\frac18,\frac38]\cup[\frac58,\frac78]$ where we can make a choice of a set $D$ which is $s$-simple, so that in the end we get $M$ to be the support of a QMF filter. Here we will make a different choice of this set $D$, and we will see that the orthonormal dilation is different from the one in Example \ref{ex1}. Here we take $D:=[\frac18,\frac38]$. Therefore 
$$M:=[0,\frac1{16}]\cup[\frac18,\frac7{16}]\cup[\frac9{16},\frac58]\cup[\frac{15}{16},1].$$
We have a partition subordinated to $M$ as follows:
$$\tau_0[0,\frac1{16}]\subset[0,\frac1{16}],\quad \tau_1[\frac{15}{16},1]\subset[\frac{15}{16},1],$$
$$\tau_1[\frac18,\frac14]\subset[\frac9{16},\frac58],\quad\tau_0[\frac14,\frac7{16}]\subset[\frac18,\frac14],\quad\tau_0[\frac9{16},\frac58]\subset[\frac14,\frac7{16}].$$
Thus, we have the following cycles $\underline0$, $\underline 1$, and $\underline{100}$, which corresponds to 
$c:=\{\theta_0=\frac17,\theta_1=\frac47,\theta_2=\frac27\}$.

The chosen paths for points in $[0,1)$ are
$$\omega(x)=\left\{\begin{array}{cc}
\underline0,&\mbox{ if }x\in[0,\frac18]\\
\underline{100},&\mbox{ if }x\in[\frac18,\frac14]\\
\underline{010},&\mbox{ if }x\in[\frac14,\frac12]\\
\underline{001},&\mbox{ if }x\in[\frac12,\frac78]\\
\underline1,&\mbox{ if }x\in[\frac78,1].
\end{array}
\right.
$$

On $[0,\frac18]$ we have $\omega(x)=\underline0$, so $d_0(\underline0)=0$, $\epsilon^{-1}([0,\frac18]\times\{\underline0\})=[0,\frac18]\times\{*\}$.

On $[\frac78,1]$ we have $\omega(x)=\underline1$, so $d_1(\underline1)=-1$, $\epsilon^{-1}([\frac78,1]\times\{\underline1\})=[\frac78,1]-1=[-\frac18,0]$.

On $[\frac18,\frac14]$ we have $\omega(x)=\underline{100}$, so $j(\underline{100})=0$, $k(\underline{100})=0$, $d_c([\frac18,\frac14]\times\{\underline{100}\})=[\frac18-\frac17,\frac14-\frac17]\times\{0\}=[-\frac1{56},\frac3{28}]\times\{0\}$.

On $[\frac14,\frac12]$ we have $\omega(x)=\underline{010}$, so $j(\underline{010})=2$, $k(\underline{010})=0$, $d_c([\frac14,\frac12]\times\{\underline{010}\})=[\frac14-\frac27,\frac12-\frac27]\times\{2\}=[-\frac1{28},\frac3{14}]\times\{2\}$.

On $[\frac12,\frac78]$ we have $\omega(x)=\underline{001}$, so $j(\underline{001})=1$, $k(\underline{001})=0$, $d_c([\frac12,\frac78]\times\{\underline{001}\})=[\frac12-\frac47,\frac78-\frac47]\times\{1\}=[-\frac1{14},\frac{17}{56}]\times\{1\}$.

Then
$$\tilde F=[-\frac18,\frac18]\times\{*\}\cup[-\frac1{56},\frac3{28}]\times\{0\}\cup[-\frac1{14},\frac{17}{56}]\times\{1\}\cup[-\frac1{28},\frac3{14}]\times 2]$$

$$\alpha(\tilde F)=[-\frac14,\frac14]\times\{*\}\cup[-\frac1{28},\frac3{14}]\times\{2\}\cup[-\frac17,\frac{17}{28}]\times\{0\}\cup[-\frac1{14},\frac37]\times\{1\}.$$

Here $\alpha(x_*,x_0,x_1,x_2)=(2x_*\mod1,2x_1\mod1,2x_2\mod1,2x_0\mod1)$. From this $\tilde P=\alpha(\tilde F)\setminus\tilde F$ and the orthogonal wavelet is
$$\tilde\psi=\chi_{\tilde P}=(\chi_{[-\frac14,-\frac18]\cup[\frac18,\frac14]},\chi_{[-\frac17,-\frac{1}{56}]\cup[\frac{3}{28},\frac{17}{28}]},\chi_{[\frac{17}{56},\frac37]},0).$$

\end{example}

\begin{example}\label{ex3}
In this example we will show that sometimes the cycles are not sufficient to describe the whole picture. Actually, for a large class of paths, we can find {\it low-pass filters for which some points will have the chosen path equal to the given path}. We will then obtain the following:

\begin{proposition}
There are low-pass filters that have chosen paths $\omega(x)$ non-eventually periodic, for a set of points $x$ of positive measure. This implies that the corresponding dilation is not realized in a sum of representations of the form $\{\hat U_C,\hat T_C\}$ with $C$ cycle.
\end{proposition}

Let $\eta=\eta_1\eta_2\dots$ be an infinite path with the property that it does not contain sequences of consecutive $0$s or consecutive $1$s of arbitrarily large lengths. (For example, any non-trivial cyclic path will have this property, or $\eta=abaabaaabaaaab\dots$ where $a=01$ and $b=10$.) Let $p-1$ be the maximum number of consecutive $0$s or $1$s that occur in $\eta$.

For a finite word $a_1\dots a_n$, let us denote by $.a_1\dots a_n:=a_1\frac12+\dots+a_n\frac{1}{2^n}$. Let $0^p$ denote a string of $p$ consecutive zeros.

Let $I$ be the interval $I:=(.10^p10,.10^p11)$.

First we claim that the intervals $(\tau_{\eta_n}\dots\tau_{\eta_1}I)_{n\geq0}$ are mutually disjoint (the interval corresponding to $n=0$ is $I$). 

Note that $\tau_{\eta_n}\dots\tau_{\eta_1}I=(.\eta_n\dots\eta_110^p10,.\eta_n\dots\eta_110^p11)$.

If $x$ is in this interval and $x$ has the binary expansion $x=a_1\frac12+a_2\frac1{2^2}+\dots$, then it is easy to see that $a_1a_2\dots$ must begin with $\eta_n\dots\eta_110^p10$. So, if $\tau_{\eta_m}\dots\tau_{\eta_1}I$ with $m>n$ intersects $\tau_{\eta_n}\dots\tau_{\eta_1}I$, then $\eta_m\dots\eta_110^p10$ must begin with $\eta_n\dots\eta_110^p10$. This implies that $\eta_{m-n}\dots\eta_110^p10$ begins with $10^p10$. But $\eta$ does not contain $p$ consecutive zeros, so the only place where we find $0^p$ is at the end. And this contradicts $m>n$.

From this it follows that the intervals $\tau_{\check\eta_n}\tau_{\eta_{n-1}}\dots\tau_{\eta_1}I$ do not intersect the intervals $\tau_{\eta_m}\dots\tau_{\eta_1}I$. (Recall that $\check\omega=1-\omega$.) Suppose by contradiction that they do intersect. Let $r(x)=2x\mod 1$ on $[0,1)$. Then 
$$\emptyset\neq r(\tau_{\check\eta_n}\tau_{\eta_{n-1}}\dots\tau_{\eta_1}I\cap\tau_{\eta_m}\dots\tau_{\eta_1}I)\subset
r(\tau_{\check\eta_n}\tau_{\eta_{n-1}}\dots\tau_{\eta_1}I)\cap r(\tau_{\eta_m}\dots\tau_{\eta_1}I)=
\tau_{\eta_{n-1}}\dots\tau_{\eta_1}I\cap\tau_{\eta_{m-1}}\dots\tau_{\eta_1}I,$$
which contradicts the previous statement.

Consider the set 
$$S:=\bigcup_{n\geq1}\tau_{\eta_n}\dots\tau_{\eta_1}I.$$

The set $S$ is $s$-simple. Indeed, if $x\in S$ then $x\in\tau_{\eta_n}\dots\tau_{\eta_1}I$ for some $n\geq 1$.
Therefore $s(x)\in\tau_{\check\eta_n}\tau_{\eta_{n-1}}\dots\tau_{\eta_1}I$ so $s(x)\notin S$.

Moreover the distance from $S$ to the boundary points $0$ and $1$ is positive. Indeed, none of the elements of the union $\tau_{\eta_n}\dots\tau_{\eta_1}I$ contains $0$ or $1$. Therefore we can consider $n$ large enough. Let $n$ be also bigger than $p$. Then the sequence $\eta_n\dots\eta_{n-p+1}$ with $n\geq p$ contains both zeros and ones
$$I_n:=\tau_{\eta_n}\dots\tau_{\eta_1}I\subset\tau_{\eta_n}\dots\tau_{\eta_{n-p+1}}[0,1)=\frac12\eta_n+\dots\frac1{2^{p}}\eta_{n-p}+\frac1{2^{p}}[0,1).$$
Since at least one of the digits $\eta_n,\dots,\eta_{n-p+1}$ is a one, it follows that $I_n\geq\frac1{2^p}$. 
Since at least one of these digits is a zero it follows that $I_m\leq \frac12+\dots+\frac1{2^{p-1}}+0+\frac1{2{p}}=1-\frac1{2^p}<1$.

Also, the distance from $S$ to $\frac12$ is positive. This is because none of the intervals contains $\frac12$ (if $\tau_{\eta_n}\dots\tau_{\eta_1}x=\frac12$ then $x=0$ or $x=1$.) And if for some $x\in I$ and some $n\geq 2$ we have
$|\tau_{\eta_n}\dots\tau_{\eta_1}x-\frac12|<\varepsilon$, then $|\tau_{\eta_{n-1}}\dots\tau_{\eta_1}x+\omega_n-1|<2\varepsilon$ so $\tau_{\eta_{n-1}}\dots\tau_{\eta_1}x$ has to be close to either $0$ or $1$, which contradicts the previous statement.

Since $S$ is $s$-simple and has positive distance to $\{0,\frac12,1\}$, we can construct an $M$ which contains $S$ and some intervals around $0$ and $1$ such that $M$ is the support of a QMF filter.

If we take $x\in I$, then $\tau_{\eta_n}\dots\tau_{\eta_1}x\in S\subset M$, so the chosen path is $\omega(x)=\eta$.

\end{example}

\begin{example}\label{ex4}
In this example, we consider the following question: in how many ways can a Parseval wavelet set be dilated? Of course, if we look at the way the low-pass filter is constructed from an MRA Parseval wavelet set, we see that there are infinitely many ways of doing that, by choosing different sets $D$ as in the proof of Proposition \ref{propm0forwavelet}.
But we would like the representations of the group $BS(1,2)$ to be different. Sure, they will have the common subrepresentation on $L^2(\br)$, but can the complementary representations be different? The answer is again yes, and we have that 

\begin{proposition}
 There are examples of MRA Parseval wavelet sets that have infinitely many orthonormal dilations with distinct representations.
\end{proposition}

Recall the notation (base two expansion): $.a_1a_2\dots a_n:=\frac12a_1+\frac1{2^2}a_2+\dots+\frac1{2^n}a_n$, and for infinite words 
$.a_1a_2\dots:=\sum_{n\geq1}\frac1{2^n}a_n$. The Euclidian order on $[0,1)$ becomes the lexicographical order on the base 2 expansions, i.e., $.a_1a_2\dots<.b_1b_2\dots$ iff for some $n\geq0$, $a_1=b_1,\dots, a_n=b_n$ and $a_{n+1}<b_{n+1}$ (there is the exception of dyadic numbers like $\frac12=.1\underline0=.0\underline1$, but these can be treated similarly).

Let $0<a<.0001$ and consider the wavelet set in Example \ref{ex1}, 
$$\hat\psi=\chi_{[-2a,-a]\cup[a,2a]}.$$

We will construct infinitely many low-pass filters $m_0=\Per(\chi_M)$ associated to this wavelet set, in such a way that their corresponding representations are distinct, and actually have only $L^2(\br)$ as the common subrepresentation.

From Example \ref{ex1}, we know that such a set $M$ must contain $\mathcal I:=[0,\frac a2]\cup[\frac12-a,\frac12-\frac a2]\cup[\frac12+\frac a2,\frac12+a]\cup[1-\frac a2,1]$ and must be disjoint from $\mathcal N:=[\frac a2,a]\cup[\frac12-\frac a2,\frac12+\frac a2]\cup[1-a,1-\frac a2]$.

We have to complete the set $\mathcal I$ with a set $D$ such that $M=\mathcal I\cup D$ gives a QMF filter, $\mathcal I\cap D=\emptyset$. We will do this in infinitely many ways, $D_n$, $n\geq 1$.

For this we consider infinitely many cycles: let $C_1$ be the cycle associated to $1001100$, $C_2$ the cycle associated to $10011001100$, $C_3$ the cycle associated to $100110011001100$ and so on. We will want $D_n$ to contain the cycle $C_n$ in its interior. 

First we have to remark a few things about the cycles $C_n$. It is easy to see that the points of the cycle $C_1$ are $.\underline{0011001}$, $.\underline{1001100}$, $.\underline{0100110}$, $.\underline{0010011}$, $.\underline{1001001}$, $.\underline{1100100}$, $.\underline{0110010}$. Note that the digits of the word associated to the cycle have to be reversed in the base $2$ expansion of the cyclic points and then cyclically permuted.

We want to make sure the cycles $C_n$ lie completely in the undecided zone 
$$\mathcal U:=[a,\frac12-a]\cup[\frac12+a,1-a].$$ Note that the point in $C_n$ closest to $0$ starts with $.001$, therefore it is bigger than $a$ (recall $a<.0001$). The point in $C_n\cap[0,\frac12]$ closest to $\frac12$ begins with $.0110$ therefore its distance to $\frac12=.01111\dots$ is at least $.0001>a$. The point in $C_n\cap[\frac12,1]$ closest to $\frac12$ begins with $.1001$ so again the distance to $\frac12=.1$ is at least $.0001>a$. Finally, the point in $C_n$ closest to $1$ begins with $.110$ so the distance to $1$ is bigger than $a$. 

Thus $C_n$ lies in the undecided zone $[a,\frac12-a]\cup[\frac12+a,1-a]$.

Next, we construct $M_n$ by adding to the set $\mathcal I:=[0,\frac a2]\cup[\frac12-a,\frac12-\frac a2]\cup[\frac12+\frac a2,\frac12+a]\cup[1-\frac a2,1]$ some intervals $I_c$, $c\in C_n$ contained in  $[a,\frac12-\frac a2]\cup[\frac12+\frac a2,1-a]$, and such that $I_c$ contains $c$ in its interior for all $c\in C_n$. 

Consider the {\it supplements} of the points in the cycle $C_n$, i.e., the points in $s(C_n)$. Since $C_n$ is in the undecided zone $\mathcal U$, and this is invariant under $s$, it follows that $C_n\cup s(C_n)$ is contained in the undecided zone. Let us call the points in $C_n\cup s(C_n)$, {\it main points}, and the points in $s(C_n)$, {\it supplements}. 

Arrange the main points on the interval. We make the following claim: the main point closest to $a$ is a supplement, the main point closest to $\frac12$ to the left of $\frac12$ is a cycle point, the main point closest to $\frac12$ to the right of $\frac12$ is a supplement, and the main point closest to $1-a$ is a supplement.

To prove the claim it is enough to prove the first and the last statement because then the other two follow by applying $s$. 

If we want a point $.a_1a_2\dots$ to be close to $0$, then we need it to start with as many zeros as possible. The cycle points start with at most two zeros: $.001$. The base two expansion for the supplements is obtained from the base two expansion of the cycle points by changing the first digit from $0$ to $1$ and vice versa. Therefore we can get the supplements to start with $.0001$ because the cycle points can start with $.1001$. Thus the one closest to $0$ will be a supplement. A similar argument works for the main point closest to $1$, and here we need the base two expansion to start with as many ones as possible.

 Let $c_l$ be the cycle point closest to $\frac12$, to the left of $\frac12$, and let $c_r$ be the cycle point closest to $\frac12$ to the right of $\frac12$. Then $s(c_l)$ will be the main point closest to $1$, and $s(c_r)$ will be the main point closest to $0$.
 
 Also, note that both intervals $[a,\frac12-a]$ and $[\frac12+a,1-a]$ contain at least $3$ cycle points (so also at least $3$ supplements). To see this we only have to count how many cycle points start with $.0$ and how many start with $.1$.

Next, consider a cycle point $c$ of $C_n$. We want to construct an interval $I_c$ associated to it, that we will add to the definition of the set $D_n$. We have two cases:

{\bf Case I.} If $c\neq c_l,c_r$, then we construct $I_c$ as follows: let $l(c),r(c)$ be the main points $l(c)<c<r(c)$ that are closest to $c$. Then $I_c:=[\frac{l(c)+c}2,\frac{c+r(c)}2]$.

{\bf Case II.} If $c=c_l$ then let $l(c_l)<c_l$ be the main point closest to $c_l$, and let $I_c=I_{c_l}:=[\frac{l(c_l)+c_l}2,\frac12-\frac a2]$. If $c=c_r$ then let $r(c_r)>c_r$ be the main point closest to $c_r$ and let $I_c:=I_{c_r}:=[\frac12+\frac a2,\frac{c_r+r(c_r)}2]$.

Note that the intervals $I_c$ are disjoint, and they are contained in $[a,\frac12-\frac a2]\cup[\frac12+\frac a2,1-a]$.
Then we define 
$$M_n:=[0,\frac a2]\cup[\frac12-a,\frac12-\frac a2]\cup[\frac12+\frac a2,\frac12+a]\cup[1-\frac a2,1]\cup\bigcup_{c\in C}I_c.$$
(Note that the union is not disjoint, because $I_{c_l}$ contains $[\frac12-a,\frac12-\frac a2]$, and $I_{c_r}$ contains $[\frac12+\frac a2,\frac12+a]$; but this will not be a problem. The set $M_n$ is indeed of the type constructed in the proof of Proposition \ref{propm0forwavelet}.)

First we have to prove that $\Per(\chi_{M_n})$ is a QMF filter, i.e., $\{M_n,s(M_n)\}$ is a partition of $[0,1)$. For this we analyze the intervals of $M_n$ in $[0,\frac12)$ and make sure that when we apply $s$ we obtain the converse situation (so that the QMF condition is satisfied).

First we have $[0,\frac a2]$ inside $M_n$ and $s([0,\frac a2])=[\frac12,\frac12+\frac a2]$ is outside $M_n$.
Then $[\frac a2,a]$ is outside $M_n$ and $s([\frac a2,a])=[\frac12+\frac a2,\frac 12+ a]$ is inside $M_n$.

Then we have the main point closest to $0$, which is $s(c_r)$. The main point closest to $c_r$ and to the right of $c_r$ is $r(c_r)$. Applying $s$ we get that the main point closest to $s(c_r)$ to the right of it is $s(r(c_r))$. The interval $[\frac12+a,\frac{c_r+r(c_r)}2]$ is in $M_n$, and its supplement $[a,\frac{s(c_r)+s(r(c_r))}2]$ is outside $M_n$ because there are no other cycle points in this region except maybe $s(r(c_r))$.

Next we consider intervals of the form $[\frac{a+b}2,\frac{b+c}2]$ where $a<b<c$ are consecutive main points in $[0,\frac12]$ and $s(c_r)<a$, $c<c_l$. If $b$ is a cycle point then this interval is contained in $M_n$. Its supplement is $[\frac{s(a)+s(b)}2,\frac{s(b)+s(c)}2]$, and $s(b)$ is a supplement. Therefore it is outside $M_n$. If $b$ is a supplement, then $s(b)$ is a cycle point and the same argument works.

Then we have the interval $[\frac{l(c_l)+c_l}2,\frac12- a]$ in $M_n$ and, using the argument that we used before for the interval $[\frac12+a,\frac{c_r+r(c_r)}2]$, its supplement is outside $M_n$.

And finally the intervals $[\frac12-a,\frac12-\frac a2]$ and $[\frac12-\frac a2,\frac12]$ can be seen to have the desired property.

This proves that $m_0=\Per(\chi_{M_n})$ is a QMF filter.

Next, we claim that for each $x\in[0,1)$ its chosen path $\omega(x)$ ends in $\underline 0$, $\underline 1$ or the infinite repetition of the finite word associated to $C_n$. It is enough to prove this for points in $M_n$, because after the first step $\tau_{\omega_1}x\in M_n$.

If $x\in[0,\frac a2]$ or $x\in[1-\frac a2,1]$ the claim is clear, the chosen path is $\underline 0$ or $\underline1$.

Suppose now $x\in I_c$ for some $c\in C_n$. Assume first that $c\neq c_l,c_r$. Since $c$ is in the cycle $C_n$, there exists $d\in C_n$ and some $i\in \{0,1\}$ such that $\tau_i c=d$. We claim that $\tau_i I_c\subset I_d$. 
We have that $I_d=[\frac{l(d)+d}2,\frac{d+r(d)}2]$ if $d\neq c_r,c_l$; if $d=c_l$ then $I_d\supset [\frac{l(c_l)+c_l}2,\frac{c_l+\frac12}2]$, because $\frac12-\frac a2>\frac{c_l+\frac12}2$; if $d=c_r$ then 
$I_d\supset [\frac{\frac12+c_r}2,\frac{c_r+r(c_r)}2]$. Therefore in all cases $I_d\supset [\frac{l(d)+d}2,\frac{d+r(d)}2]$ where $r(c_l):=\frac12=:l(c_r)$. Note also that $I_d$ is completely contained in $[0,\frac12]$ or $[\frac12,1]$ so $R(x)=2x\mod1$ is injective on $I_d$. Since $\tau_i c=d$, we have that $R(x)=2x-i$ on $I_d$, so the inverse of $R$ on $I_d$ is $\tau_i$.

Note that if $x$ is a main point then $R(x)$ is a cycle point. This is clear for cycle points; for supplements, $s(x)$ is a cycle point, and $R(x)=R(s(x))$. So if $d\neq c_r$, then $R(l(d))$ is a cycle. If $d=c_r$ then $R(l(c_r))=R(\frac12)=0$. A similar argument works for $R( r(d))$ which is a cycle point or $1$. Then $R(\frac{l(d)+d}2)=\frac{R(l(d))+R(d)}2=\frac{R(l(d))+c}2\leq\frac{l(c)+c}2$, because the first is the midpoint between two cycle points (or perhaps $0$ and a cycle point), and the last is the midpoint between two main points. Similarly for $R(\frac{d+r(d)}2)\geq\frac{c+r(c)}2$. This shows that $R(I_d)\supset I_c$. Taking the inverse we obtain $\tau_i(I_c)\subset I_d$. Note that when $d=c_l$, we actually have $\tau_i(I_c)\subset\tilde I_{c_l}:=[\frac{l(c_l)+c_l}2,\frac{c_l+\frac12}2]$, and when $d=c_r$, $\tau_i(I_c)\subset\tilde I_{c_r}:=[\frac{\frac12+c_r}2,\frac{c_r+r(c_r)}2]$.

Consider now $I_{c_l}$. Since $c_l$ is on the cycle $C_n$, there is some $i\in\{0,1\}$ and some $d\in C_n$ such that $\tau_ic_l=d$. Note that $d$ cannot be $c_l$ or $c_r$. It is not $c_l$ because $c_l$ is not a fixed point, the cycle $C_n$ is longer than 1. It is not $c_r$ because $c_l$ starts with $.0110$ so $\tau_ic_l$ starts with $.00110$ or $.10110$. But $c_r$ starts with $.1001$. We claim that, with $\tilde I_{c_l}:=[\frac{l(c_l)+c_l}2,\frac{c_l+\frac12}2]$, we have $\tau_i\tilde I_{c_l}\subset I_d$. As before, we look at $R$ applied to the endpoints of $I_d$, and we have already seen that $R(\frac{l(d)+d}2)\leq\frac{l(c_l)+c_l}2$. As before
$R(\frac{d+r(d)}2)=\frac{R(d)+R(r(d))}2=\frac{c_l+R(r(d))}2$ and $R(r(d))$ is a cycle point to the right of $R(d)=c_l$. Thus $R(r(d))>\frac12$. This shows $R(I_d)\supset \tilde I_{c_l}$, so $\tau_i(\tilde I_{c_l})\subset I_d$.

A similar argument works for $c_r$: if $\tau_ic_r=d$ for some $i\in \{0,1\}$ and $d\in C_n$, and if we denote by 
$\tilde I_{c_r}:=[\frac{\frac12+c_r}2,\frac{c_r+r(c_r)}2]$ then $\tau_i(\tilde I_{c_r})\subset I_d$.

Let $\tilde I_c:=I_c$ for $c\neq c_l,c_r$. We have that, if $\tau_ic=d$, then for all $x\in \tilde I_c$, one has $\tau_ix\in \tilde I_d\subset M_n$. This means that the chosen path of $x$ starts with $i$, and by induction, the digits of the chosen path will coincide with the digits that determine the cycle $C_n$.

The only remaining intervals are $J_l:=[\frac{c_l+\frac12}2,\frac12-\frac a2]$ and $J_r:=[\frac12+\frac a2,\frac{\frac12+c_r}2]$. If $x\in J_l$ then $c_l<x<\frac12$, therefore $x$ starts with $.011$. If $x\in J_r$ then $\frac12<x<c_r$. Hence $x$ starts with $.100$. Thus
if $x\in J_l$ then $\tau_ix$ starts with $.0011$ or $.1011$ so it cannot be in $J_l$ or $J_r$. Therefore, if $x\in J_l\cup J_r$ and $\omega(x)=\omega_1\omega_2\dots$ is its chosen path, then $\tau_{\omega_1}x$ is in $M_n$ but cannot be in $J_l$ or $J_r$. This implies that $\tau_{\omega_1}x$ fits into one of the previous cases, so its chosen path ends in a repetition of one of the cycles. A similar argument works for $J_r$.

Hence, every chosen path will end in $\underline0$, $\underline1$ or the infinite repetition of the word associated to the cycle $C_n$. It is also clear that $C_n$ lies in the interior of $M_n$ and the conditions of Proposition \ref{propchardense} are satisfied. Then Theorem \ref{thdilmra} and the proof of Theorem \ref{th3_24}(ii) shows that the orthonormal dilation of our Parseval wavelet set has the representation $\{U_0\oplus U_{C_n}, T_0\oplus T_{C_n}\}$. With Lemma \ref{lem}, we see that the representations $\{U_{C_n},T_{C_n}\}$ are mutually disjoint. This proves our claims.

\end{example}

\begin{acknowledgements}
We would like to thank Professors Uffe Haagerup and Palle Jorgensen for their suggestions and ideas. We also thank the anonymous referee for his comments and suggestions which helped improve our paper. 
\end{acknowledgements}

\bibliographystyle{alpha}
\bibliography{dilws}
\end{document}